\documentclass{article}

\usepackage[utf8]{inputenc}

\usepackage{amsthm,amssymb,stmaryrd}

\usepackage{graphicx}

\usepackage[cmtip,all]{xy}

\usepackage{fdsymbol}
\usepackage{dsfont}
\usepackage{mathrsfs}
\usepackage{mathtools}
\usepackage{tikz}

\usetikzlibrary{decorations.pathmorphing,decorations.pathreplacing}
\usepackage[outline]{contour}
\contourlength{2pt}

\usetikzlibrary{patterns}

\usepackage{hyperref}
\usepackage{epigraph}
\usepackage{accents}
\usepackage{xfrac}
\usepackage{enumitem}
\usepackage{stackrel}

\newtheorem{theorem}{Theorem}[section]

\theoremstyle{definition}
\newtheorem{definition}[theorem]{Definition}
\newtheorem{example}[theorem]{Example}

\theoremstyle{remark}

\numberwithin{equation}{section}

\title{Recognition of commutative algebra spectra through an idempotent quasiadjunction}
\author{Renato Vasconcellos Vieira\thanks{This work was financed by the grant 2020/06159-5, São Paulo Research Foundation (FAPESP)}}
\date{January 2021}

\begin{document}

\maketitle

\begin{abstract}
In this article a recognition principle for $\infty$-loop pairs of spaces of connective commutative algebra spectra over connective commutative ring spectra is proved. This is done by generalizing the classical recognition principle for connective commutative ring spectra using relative operads. The machinery of idempotent quasiadjunctions is used to handle the model theoretical aspects of the proof. 
\end{abstract}

\section{Introduction}

The category $\texttt{Sp}_{\mathds N}$ of sequential prespectra \cite{Lima58} consists of sequences of spaces $\langle Y_N\rangle\in\prod_{\mathds N}\texttt{Top}_\ast$ equipped with structural maps $\sigma^M_N:Y_M\wedge\mathds S^{N-M}\rightarrow Y_N$ for $M\leq N$ satisfying compatibility conditions. An $\Omega$-spectrum is a prespectrum whose adjoint structural maps $\tilde\sigma^M_N$ are weak equivalences, which by Brown representability represents (co)homology theories \cite{Brown62}. Spectra are prespectra such that the $\tilde \sigma^N_M$ are homeomorphisms (see for instance \cite{EKMM}). In this article we will work exclusively in the category of prespectra, so from now on we will simply refer to prespectra as spectra. From $\texttt{Sp}_{\mathds N}$ we can define via filtered colimits over the dual structural maps the $\infty$-loop spaces functor
$$
\Omega^\infty:\texttt{Sp}_{\mathds N}\rightarrow \texttt{Top}_\ast;\hspace{0.5cm}\Omega^\infty Y:=\text{colim}_{\mathds N} Y_{N}^{\mathds S^{N}}.
$$

The $\infty$-loop spaces $\Omega^\infty Y$ are homotopy commutative $H$-spaces, but such description ignores a lot of information. In order to describe the algebraic structure completely we require an $E_\infty$-operad $\mathcal E$, a gadget used to describe topological spaces with operations that are associative and commutative up to coherent homotopy \cite{MayGeomItLoopSp}.

For $\mathbb S$ the category of finite sets a topological operad is a contravariant functor equipped with an abstract identity element and composition maps 
\begin{gather*}
    \mathcal P:\mathbb S^{\text{op}}\rightarrow \texttt{Top};\\
    \textbf{id}\in \mathcal P\underline 1,
    \hspace{1cm}\circ:\mathcal P A\times \textstyle\prod_{A}\mathcal PB^a\rightarrow \mathcal P\textstyle\Sigma_{A} B^a
\end{gather*}
with $\mathcal P\emptyset=\ast$ satisfying invariance, associativity and unitary laws. We can interpret points in the underlying spaces as abstract multivariable functions with inputs indexed by the sets $A$. This structures allows us to define via the coend construction \cite{Lo15} the monad 
\begin{gather*}
    P:\texttt{Top}_\ast\rightarrow \texttt{Top}_\ast, \hspace{1cm} PX:=\textstyle\int^{\mathbb S^{\text{inj}}}\mathcal PA\times X^{\times A};\\
    \eta x:=[\textbf{id},x],
    \hspace{1cm} \mu [\alpha,\langle[\beta^a,\langle x^{a,b}\rangle]\rangle]:=[\alpha\circ\langle \beta^a\rangle,\langle x^{a,b}\rangle].
\end{gather*}

The category $\mathcal P[\texttt{Top}]$ of $\mathcal P$-spaces consists of pointed spaces $X\in\texttt{Top}_\ast$ equipped with maps $\xi:PX\rightarrow X$ compatible with the monad maps, which we interpret as an instantiation of the abstract operations of $\mathcal P$.

An important family of operads are the embeddings operads $\text{Emb}_N$ for $N\in\mathds N$ with 
$$
    \text{Emb}_NA:=\left\{\alpha=\langle \alpha_a\rangle\in (\mathds R^N)^{\sqcup_A\mathds R^N}\mid \alpha\text{ an embedding}\right\}.
$$
There are natural operad inclusions $\text{End}_M\hookrightarrow \text{End}_N$ and we define $\text{End}_\infty:=\text{colim}_{\mathds N}\text{End}_N$. All $N$-loop spaces are naturally $\text{End}_N$-spaces with
$$
    \alpha\langle \gamma^a\rangle:=\left(\vec u\mapsto \begin{cases}
        \gamma^a\alpha_a^{-1}\vec u,&\vec u\in\alpha_a\mathds R^N\\
        \ast,&\vec u\not\in\alpha\sqcup_A\mathds R^N
    \end{cases}\right)
$$
and these induce $\text{End}_\infty$-space structures on $\infty$-loop spaces.

An $E_\infty$-operad is an operad $\mathcal E$ with each underlying space $\mathcal E A$ a contractible free $\mathbb SA$-space. For the purpose of studying $\infty$-loop spaces we further require $E_\infty$-operads to be equipped with an operad map $\psi:\mathcal E\rightarrow \text{End}_\infty$. This allows us to define by pullback an $\mathcal E$-algebra structure on $\infty$-loop spaces $\Omega^\infty Y$, which induces a functor $\Omega^\infty:\texttt{Sp}\rightarrow \mathcal E[\texttt{Top}]$. This functor is not a right adjoint since any abelian group $G$ is an $\mathcal E$-spaces and it can be shown that due to the strictness of the operations in $G$ any $\mathcal E$-map $\varphi:G\rightarrow \Omega^\infty Y$ must be trivial, so no unit of adjunction can be constructed.

In May's recognition theorem \cite{MayGeomItLoopSp} the solution was to consider the resolution of $\mathcal E$-spaces by the bar construction
$$
    \overline B:\mathcal E[\texttt{Top}]\rightarrow \mathcal E[\texttt{Top}];
    \hspace{1cm}\overline B X:=B(E,E,X),
$$
which comes equipped with a natural weak equivalence $\eta':\overline B\Rightarrow Id$.

The maps $\psi:\mathcal E\rightarrow \text{End}_\infty$ induce by pullback a suboperad filtration $\mathcal E_N$ on $\mathcal E$. If each underlying space $\mathcal E_NA$ is equivariantly homotopy equivalent to the configuration space of $A$ elements in $\mathds R^N$ then we can define the $\infty$-delooping functor
$$
    B^\infty:\mathcal E[\texttt{Top}]\rightarrow \texttt{Sp};
    \hspace{1cm}B^\infty X:=\langle B(\Sigma^N,E_N,X)\rangle;
$$
such that there is a natural transformation $\eta:\overline B\Rightarrow \Omega^\infty B^\infty$, with $\eta_X$ a weak equivalence if and only if $X$ is grouplike, meaning that $\pi_0 X$ is not only a monoid but also a group.

Dually there is no counit map. There is a spectrification functor\footnote{The spectrification functor $\widetilde \Omega$ is the left adjoint to the inclusion of the category of spectra in the sense used in \cite{EKMM} into the category of prespectra, hence the name.}
$$
    \widetilde\Omega:\texttt{Sp}\rightarrow\texttt{Sp};\hspace{1cm}\widetilde\Omega Y:=\langle\text{colim}_{M\leq N}\widetilde Y_N^{\mathds S^{N-M}}\rangle,
$$
where $\widetilde Y$ is a certain inclusion prespectrum constructed from $Y$, such that we have natural inclusions $\epsilon':Id\Rightarrow\widetilde\Omega$ which are stable weak equivalences. This functor plays an important role in the construction of the stable model structures of spectra. There is a natural transformation $\epsilon: B^\infty\Omega^\infty\Rightarrow\widetilde\Omega$ such that the equation $\Omega^\infty\epsilon \eta_{\Omega^\infty}= \Omega^\infty\epsilon'\eta'_{\Omega^\infty}$ holds in $\mathcal E[\texttt{Top}]$ and we have a homotopy equivalence $\epsilon_{B^\infty}B^\infty\eta_X\simeq \epsilon'_{B^\infty X}B^\infty\eta_X'$ in $\texttt{Sp}$.
\begin{equation*}
    \xymatrix@C=1cm@R=0.4cm{
        B^\infty\overline BX\ar[d]_{B^\infty\eta'_X}^\sim\ar[r]^{B^\infty\eta_X}&B^\infty\Omega^\infty B^\infty X\ar[d]^{\epsilon_{B^\infty X}}\\
        B^\infty X\ar[r]^\sim_{\epsilon'_{B^\infty X}}&\widetilde\Omega B^\infty X\\
    }\hspace{1cm}\xymatrix@C=1cm@R=0.4cm{
        \overline B \Omega^\infty Y\ar[d]^\sim_{\eta'_{\Omega^\infty Y}}\ar[r]^{\eta_{\Omega^\infty Y}}&\Omega^\infty B^\infty\Omega^\infty Y\ar[d]^{\Omega^\infty \epsilon_Y}\\
        \Omega^\infty Y\ar[r]^\sim_{\Omega^\infty\epsilon'_Y}&\Omega^\infty\widetilde\Omega Y
    }
\end{equation*}

Note the similarity of these equations to the ones for an adjunction. Indeed if $\overline B$, $\widetilde \Omega$, $\eta'$ and $\epsilon'$ were substituted by identities and both equations held strictly we would have an adjunction in the regular sense. In \cite{Vi20} I defined a generalization of Quillen adjunctions, called weak Quillen quasiadjunctions, that allowed for units and counits to exist up to functorial resolutions. I proved that weak Quillen quasiadjunctions still induce adjunctions of the homotopy categories, generalizing the analogous result for Quillen adjunctions. In the same vein I defined a generalization of Quillen idempotent (co)monads that induce left (right) Bousfield localizations of model structures, and through these we have a natural definition of idempotent quasiadjunctions which induce equivalences between the associated homotopy subcategories.

Adapting May's original proof of the recognition principle in \cite{MayGeomItLoopSp,MayEspcGrpComplPermCat} we can show that the weak Quillen quasiadjunction
$$
    (B^\infty\dashv_{\ \overline B,\widetilde\Omega}\Omega^\infty):\mathcal E[\texttt{Top}]\leftrightharpoons \texttt{Sp}_{\mathds N}
$$
is idempotent and induces an equivalence between the homotopy category of grouplike $\mathcal E$-spaces and the the homotopy category of connective spectra. Idempotent quasiadjunctions provide a model categorical axiomatization of the essential elements of May's original proof, and it can be adapted to prove variations of the recognition principle. For instance the relative recognition principle for $\infty$-loop pairs of spaces of spectra maps of degree 1 was proved using the above machinery and relative operads in \cite{Vi20}.

In this article I show that the machinery of quasiadjunctions and relative operads are also compatible with actions by a natural relative version $\mathscr L^\shortrightarrow$ of the linear isometries operad $\mathscr L$, introducing in particular a relative version of actions between operads which provides a natural definition of $E_\infty$-algebra spaces over $E_\infty$-ring spaces. The main theorem \ref{RecogCAlgSp} is a recognition principle for $\infty$-loop pairs of spaces of commutative algebra spectra over commutative ring spectra. Explicitly it states that the homotopy category of algebralike\footnote{An $E^\shortrightarrow_\infty$-algebra is algebralike if it admits additive inverses up to homotopy.} $E^\shortrightarrow_\infty$-algebras is equivalent to the homotopy category of connective commutative algebra spectra over connective commutative ring spectra. As in \cite{EKMM} this will require us to work on the more structured category $\texttt{Mod}_{\mathds S}$ of $\mathds S$-modules, which is monoidal and so provides a convenient language to describe algebraic structures. In particular we will work with the coordinate-free spectra of \cite{LMS} which substitutes the natural numbers $\mathds N$ by the set of finitely dimensional subspaces of some countably infinite dimensional inner product spaces such as $\mathds R^\infty$ as the indexing set of spectra. This result is a simple consequence of the intermediary theorems \ref{MorphRecogPrin} and \ref{Idempotent}, which is a recognition principle for $\infty$-loop pairs of spaces of spectra maps.

\subsection{Structure of the article}

In section 2 we review the definition of weak Quillen quasiadjunctions, idempotent quasimonads and idempotent quasiadjunctions. Our main theorem will be a particular case of the fact that idempotent quasiajunctions induce equivalences between the associated homotopy subcategories.

In section 3 we present the definition of $E_\infty$-algebras over $E_\infty$-rings through relative operads. A detailed description of relative sets and filtered rooted relative trees and operations on them will be required to construct bar resolutions and delooping spectra, as well as describe their algebraic structures. We then give a brief review on relative operads and $E^\shortrightarrow_\infty$-operads and the bar resolution of $E^\shortrightarrow_\infty$-algebras. Relative operad actions are then introduced which provides an account of distributivity laws between multiplicative and additive relative operad actions and is central in the definition of the category $(\mathscr E^\shortrightarrow,\mathscr L^\shortrightarrow)[\texttt{Top}]$ of $E_\infty$-algebras. We also give a brief review of how the Quillen model structure on $(\mathscr E^\shortrightarrow,\mathscr L^\shortrightarrow)[\texttt{Top}]$ is transfered from the one on $\texttt{Top}^2_\ast$.

The main theorems are in section 4. We review the basics of coordinate-free spectra and the construction of the stable mixed model structure. The recognition principle for $\infty$-loop pairs of spaces of spectra maps is proved via an idempotent quasiadjunction in theorems \ref{MorphRecogPrin} and \ref{Idempotent}, which imply the homotopy category of grouplike $E^\shortrightarrow_\infty$-pairs is equivalent to the homotopy category of spectra maps between connective spectra. After a review of the basics of $\mathds S$-modules and commutative algebra spectra, including the construction of stable mixed model structures, the main theorem \ref{RecogCAlgSp} is proved.

\subsection{Notation and terminology}

We assume the theory of model categories in \cite{GJ09, Hi09, Ho07}, and the theory of monoids, their algebras and the bar construction in \cite[Section 9]{MayGeomItLoopSp}. In diagrams in a model category $\mathcal T$ the morphisms in the class of weak equivalences $W$ are denoted by arrows marked with a tilde $\xrightarrow{\sim}$, the ones in the class of cofibrations $C$ by hooked arrows $\hookrightarrow$ and the ones in the class of fibrations $F$ by double headed arrows $\twoheadrightarrow$. The functorial weak factorization systems are denoted by $(\text{Fat}_{C,F_t},C-,F_t-)$ and $(\text{Fat}_{C_t, F},C_t-,F-)$ such that a morphism $f\in\mathcal T(X,Y)$ is factored for instance as $\xymatrix@C=0.6cm{X\ar@{^{(}->}[r]^(0.37){Cf}&\text{Fat}_{C,F_t} f \ar@{->>}[r]^(0.63){F_tf}_(0.63)\sim&Y}$.

The notations $\mathfrak{C}:\mathcal T\rightarrow\mathcal T$ and $\text{cof}:\mathfrak{C}\Rightarrow Id$ are used for the cofibrant resolution functor and the associated natural trivial fibration, and the notations $\mathfrak F:\mathcal T\rightarrow \mathcal T$ and $\text{fib}:Id\Rightarrow \mathfrak{F}$ are used for the fibrant resolution functor and the associated natural trivial cofibration. The homotopy category of $\mathcal T$ with objects the bifibrant objects of $\mathcal T$ and with morphisms between bifibrant objects $X$ and $Y$ the set $\mathcal T(X,Y)/_\simeq$ of homotopy classes of maps \cite[Section 1.2]{Ho07} is denoted by $\mathcal Ho\mathcal T$. 

The monoidal category $\texttt{Top}$ of compactly generated weakly Hausdorff spaces as presented in Strickland's \cite{CGWHspaces} admits two model structures, the cofibrantly generated Quillen model structure \cite{QuHomAlg} with weak equivalences the weak homotopy equivalences ($q$-equivalences), fibrations of this model structure the Serre fibrations ($q$-fibrations), and cofibrations retracts of inclusions of well pointed relative CW-complexes ($q$-cofibrations), and the Hurewicz/Str\o m model structure \cite{Strom72} with  distinguished classes of maps the homotopy equivalences ($h$-equivalences), the Hurewicz fibrations ($h$-fibrations) and the Hurewicz cofibrations ($h$-cofibrations). As Cole proved in \cite{Cole06} we can mix these model structures into one with weak equivalences the $q$-equivalences, fibrations the $h$-fibrations and cofibrations the maps that can be factored as a $q$-cofibration followed by an $h$-equivalence. We use the notation $K\subset_{\text{cpct}} X$ to indicate $K$ is a compact subspace of $X$. We denote by $I$ the interval $[0,1]\subset\mathbb R$. We denote by $\texttt{Top}_\ast$ the category of pointed spaces and for $X\in \texttt{Top}$ we denote by $X_+\in \texttt{Top}_\ast$ the pointed space obtained by adjoining a disjoint base point.

We denote by $\mathcal T^\shortrightarrow$ the category of morphisms $f:X_d\rightarrow X_c$ in $\mathcal T$ as objects and commutative squares as morphisms. For notational convenience we denote elements of categories of pairs $\mathcal T^2$ as $X=(X_d,X_c)$, and we will consider relative operads colored on the set $\{d,c\}$, with $d$ being the ``domain'' color and $c$ the ``codomain'' color.

Let $\mathscr I$ denote the topological category of finite or countably infinite dimensional real inner product spaces and linear isometries, with the topology defined as the colimit of the finite dimensional sub-spaces. This category is monoidal under direct sums. For $\mathds U\in\mathscr I$ we denote by $\mathscr A_{\mathds U}$ the set of finite dimensional subspaces of $\mathds U$, partially ordered by inclusion, and for $U\in\mathscr A_{\mathds U}$ we define $\mathscr A_U:=\{V\in\mathscr A_{\mathds U}\mid U\leq V\}$. For $\mathds U=\mathds R^\infty$ we simply write $\mathscr A:=\mathscr A_{\mathds R^\infty}$. For $\langle f_a\rangle\in\mathscr I(\oplus_A\mathds U^a,\mathds V)$ and $\langle \vec u^a\rangle \in \oplus_A\mathds U^a$ we use the Einstein summation convention $f_a\vec u^a:=\sum_Af_a\vec u^a$. For $\mathds U\in\mathscr I$ and $U\subset \mathds U$ any subspace we use the notation $U^\bot:=\{\vec v\in\mathds U\mid \forall \vec u\in U:\vec v\cdot\vec u=0\}$ for the orthogonal complement. For any $\vec u\in \mathds U$ and $U\subset \mathds U$ we use the notation $\vec u_U$ to denote the projection of $\vec u$ on $U$. For $\vec v\in \mathds V$ and $f\in\mathscr I(\mathds U,\mathds V)$ we use the notation $\vec v_f:=\vec v_{f\mathds U}\in\mathds V$ for the projection of $\vec v$ onto the image of $f$ and $\vec u
^f:=f^{-1}\vec v_f\in\mathds U$. For all $U\in \mathscr A_{\mathds U}$ let $\mathds S^U$ be the one point compactification of $U$ obtained by adding a point $\infty$ at infinity and for $(U,V)\in \Sigma_{\mathscr A}\mathscr A_U$ let $V-U:=U\cap V^\bot$.

We will make extensive use of mapping spaces $Y^X$ and will express their elements as $x\mapsto \Phi$ for some expression $\Phi$ which may use the variable $x$. For $X$ a set (or space) equipped with an equivalence relation $\sim$ we will denote the equivalence classes of $x\in X$ using square brackets $[x]\in X/_\sim$.

We denote by $\texttt{Set}$ the category of sets and functions, by $\mathbb S^\text{inj}$ the subcategory of finite sets and injections and by $\mathbb S$ the subcategory of finite sets and bijections. We will use the notation $\underline m$ for the sets $\{1,\dots,m\}$, with $\underline 0=\emptyset$.

Given a class $A$ and a family of classes $\langle B^a\rangle$ indexed by $A$ the dependent sum $\Sigma_AB^a$ is the class of pairs $(a,b)$ with $a\in A$ and $b\in B^a$ and the dependent product $\Pi_AB^a$ is the class of sequences $\langle b^a\rangle$ indexed on $A$ with $b^a\in B^a$ for each $a\in A$, or equivalently it is the class of sections of the natural surjection $\Sigma_AB^a\rightarrow A$.

We denote by $\texttt{POSet}$ the category of ordered sets and monotone functions and $\Delta$ the full subcategory on $\langle m\rangle=\langle 0<\cdots<m\rangle$ for $m\in\mathds N$. This category is generated by the coface injections $\partial_i:\langle m-1\rangle\rightarrow\langle m\rangle$, with $i\not\in \partial_i\langle m-1\rangle$, and codegeneracy surjections $\delta_i:\langle m+1\rangle\rightarrow\langle m\rangle$, with $\delta_ii=\delta_i(i+1)$, for all $i\in\langle m\rangle$.

Consider the cosimplicial space of partitions of the interval $\text{Part}^-\in\texttt{Top}^{\Delta}$ with
\begin{gather*}
    \text{Part}^{\langle m \rangle}:=\texttt{POSet}(\langle m-1\rangle,I);\\
    \partial_i\cdot t:=\begin{cases}
        \left\langle\begin{cases}
        0,&j=0\\
        t^{j-1},&j> 0
    \end{cases}\right\rangle,&i=0\\
    \left\langle\begin{cases}
        t^j,&j< i\\
        t^{j-1},&j\leq i
    \end{cases}\right\rangle,&0<i<m\\
    \left\langle\begin{cases}
        t^j,&j<m-1\\
        1,&j=m-1
    \end{cases}\right\rangle,&i=m
    \end{cases},
    \hspace{1cm}
    \delta_i\cdot t:=\left\langle\begin{cases}
        t^j,&j< i\\
        t^{j+1},&j\geq i
    \end{cases}\right\rangle,
\end{gather*}
with $\text{Part}^{\langle m \rangle}$ topologized as a subspace of $I^{\langle m-1\rangle}$. For each $\langle t^{a}\rangle\in\Pi_A\text{Part}^{\langle m ^a\rangle}$ there is a unique
$$
    \textstyle\lhd_A t^a\in \text{Part}^{\langle\sum_A m ^a\rangle}
$$
obtained by ordering the elements $t^{a,i}$ for $a\in A$ and $i\in\langle m ^a-1\rangle$. For each $a\in A$ and $\langle  t^{a'}\rangle\in \Pi_A\Delta^{\langle m ^{a'}\rangle}$ we can define
$$
    \textstyle \delta^a\in \Delta(\langle\sum_A m^{a'}\rangle,\langle m^a\rangle),
    \hspace{0.15cm}
    \delta^ai:=\begin{cases}
        \min(j\mid (\lhd_A t^{a'})^i \leq t^{a,j}),
        &(\lhd_A t^{a'})^i\leq t^{a,m^a-1}\\
        m^a,
        &(\lhd_A t^{a'})^i>t^{a,m^a-1}
    \end{cases}
$$
such that $\delta^a\cdot\lhd_A t^{a'}= t^a$.

For any simplicial space $X^-\in\texttt{Top}^{\Delta^{\text{op}}}$ its geometric realization $|X^-|$ is defined via the coend construction \cite{Lo15} as
$$
    \textstyle |X^-|:=\int^\Delta X^{\langle m\rangle}\times \text{Part}^{\langle m\rangle}.
$$

The reason we consider the geometric realization via the partitions cosimplicial space instead of the usual homeomorphic cosimplicial space of topological simplexes is that this choice simplifies the algorithm in \cite[Theorem 11.5]{MayGeomItLoopSp}.

\section{Idempotent quasiadjunctions}

\subsection{Weak Quillen quasiadjunction}

The following definition introduced in \cite{Vi20} is a generalization of Quillen adjunctions between model categories. The basic idea is that to construct the unit and counit natural transformations of an adjuction between the homotopy categories it suffices to construct a unit natural span and counit natural cospan at the model categories level, plus some natural compatibility conditions with the model structures.

\begin{definition}\label{QsiAdjQuilFr} Let $\mathcal T$ and $\mathcal A$ be model categories. A \textit{weak Quillen quasiadjunction}, or just \textit{quasiadjunction}, between $\mathcal T$ and $\mathcal A$, denoted by 
$$
(S \dashv_{\ \mathscr{C},\mathscr{F}}\Lambda):\mathcal T\rightleftharpoons \mathcal A,
$$
is a quadruple of functors
$$
\xymatrix{
    \mathcal T\ar@<0.1cm>[r]^{S}\ar@(dl,ul)[]^{\mathscr C}&\mathcal A\ar@<0.1cm>[l]^{\Lambda}\ar@(ur,dr)[]^{\mathscr F}
}
$$
with $S$ the \textit{left quasiadjoint} and $\Lambda$ the \textit{right quasiadjoint}, equipped with a \allowdisplaybreaks natural span in $\mathcal T$ and a natural cospan in $\mathcal A$
$$
\xymatrix@=0.5cm{
Id_{\mathcal T}
&\mathscr C\ar@{=>}[l]_(0.4){\eta'}^(0.4)\sim\ar@{=>}[r]^(0.4)\eta
&\Lambda S
&
&S\Lambda\ar@{=>}[r]^\epsilon
&\mathscr F
&Id_{\mathcal A}\ar@{=>}[l]_{\epsilon'}^\sim
}$$
such that
\begin{enumerate}[label=(\roman*)]
    \item $S$ is left derivable;
        
    \item $\Lambda$ is right derivable;
        
    \item $\mathscr{C}$ and $\mathscr{F}$ preserve cofibrant and fibrant objects;
        
    \item $\eta'$ and $\epsilon'$ are natural weak equivalences;
        
    \item If $X\in\mathcal T$ is cofibrant then $\epsilon_{SX}S\eta_X\simeq \epsilon'_{SX}S\eta_X'$;
        
    \item If $Y\in\mathcal A$ is fibrant then $\Lambda\epsilon_Y \eta_{\Lambda Y}\simeq \Lambda\epsilon'_Y\eta'_{\Lambda Y}$.
\end{enumerate}
\begin{equation*}
        \xymatrix@C=0.8cm@R=0.4cm{
            S\mathscr CX\ar[d]_{S\eta'_X}^\sim\ar[r]^{S\eta_X}&S\Lambda S X\ar[d]^{\epsilon_{SX}}\\
            SX\ar[r]^\sim_{\epsilon'_{SX}}&\mathscr F SX
        }\hspace{1cm}\xymatrix@C=0.8cm@R=0.4cm{
            \mathscr C \Lambda Y\ar[d]^\sim_{\eta'_{\Lambda Y}}\ar[r]^{\eta_{\Lambda Y}}&\Lambda S\Lambda Y\ar[d]^{\Lambda \epsilon_Y}\\
            \Lambda Y\ar[r]^\sim_{\Lambda\epsilon'_Y}&\Lambda\mathscr F Y\\
        }
    \end{equation*}
\end{definition}

\begin{theorem}[{\cite[Theorem 2.1.2]{Vi20}}]
    A quasiadjunction induces an adjunction
\begin{gather*}
    (\mathbb LS \dashv\mathbb R\Lambda):\mathcal Ho \mathcal T\rightleftharpoons \mathcal Ho \mathcal A;\\
    \xymatrix@C=1.8cm{
    Id_{\mathcal Ho\mathcal T}\ar@{=>}[r]^{[\text{cof}\eta'_{\mathfrak C}]^{-1}}
    &\mathbb L\mathscr C
    \ar@{=>}[r]^(0.45){[(\Lambda\text{fib}_S\eta)_{\mathfrak C }]}
    &\mathbb R\Lambda\mathbb LS},
    \\
    \xymatrix@C=1.8cm{\mathbb LS\mathbb R\Lambda
    \ar@{=>}[r]^{[(\epsilon S\text{cof}_\Lambda)_{\mathfrak F }]}
    &\mathbb R\mathscr F\ar@{=>}[r]^{[\epsilon'_{\mathfrak F}\text{fib}]^{-1}}
    &Id_{\mathcal Ho\mathcal A}
}
\end{gather*}
between the homotopy categories.
\end{theorem} 

\subsection{Idempotent quasi(co)monads}

The following generalization of idempotent Quillen monads \cite{BF78} was also introduced following the same principle of only requiring the existence of a unit natural span, and they also induce Bousfield localizations.

\begin{definition}\label{QsiMonIdempQuil}
    Let $\mathcal T$ be a right proper model category with distinguished subclassses of morphisms $(W,C,F)$. A \textit{Quillen idempotent quasimonad on $\mathcal T$}, or simply an \textit{idempotent quasimonad}, is a pair of endofunctors $Q,\overline{\mathscr C}:\mathcal T\rightarrow\mathcal T$ equipped with a natural span 
    $$
        \xymatrix@=0.5cm{
    Id_\mathcal T&\overline{\mathscr C}\ar@{=>}[l]_(0.4){\eta'}^(0.4)\sim\ar@{=>}[r]^(0.4)\eta&Q
    }
    $$
    such that:
    \begin{enumerate}[label=(\roman*)]
        \item $Q$ preserves weak equivalences;
        
        \item $Q\eta$ and $\eta_{Q}$ are natural weak equivalences;
        
        \item If $f\in\mathcal T(X,B)$, $p\in F(E,B)$ and $\eta_E,\eta_B,Qf\in W$ then $Q(f^\ast p)\in W$;
        $$\xymatrix@R=0.5cm{
                X\times_BE\ar[d]_{p^\ast f}\ar[r]^(0.6){f^\ast p}&E\ar@{->>}[d]_p&\overline{\mathscr C}E\ar[d]_{\overline{\mathscr C}p}\ar[r]^{\eta_E}_\sim\ar[l]_{\eta'_E}^\sim&
                QE\ar[d]_{Qp}&Q(X\times_BE)\ar[d]^{Q (p^\ast f)}\ar[l]_(0.6){Q(f^\ast p)}^(0.6)\sim\\
                X\ar[r]_f&B&\overline{\mathscr C}B\ar[l]^{\eta'_B}_\sim\ar[r]_{\eta_B}^\sim&QB&QX\ar[l]^{Qf}_\sim
            }$$
            
            \item $\eta'$ is a natural weak equivalence;
    
            \item If $\iota\in C(\overline{\mathscr C} X,K)$ then $\iota_\ast \eta'\in W$.
        $$\xymatrix@R=0.5cm{
            \overline{\mathscr C}X\ar@{^{(}->}[d]_\iota\ar[r]^{\eta'}_\sim&X\ar[d]^{\eta'_\ast\iota}\\
            K\ar[r]^(0.35)\sim_(0.35){\iota_\ast\eta'}&K\sqcup_{\overline{\mathscr C}X}X
        }$$
    \end{enumerate}
\end{definition}

\begin{theorem}[{\cite[Theorems 2.3.5 and 2.3.6]{Vi20}}]
    An idempotent quasimonad induces a left Bousfield localization
$$
    \mathcal T_Q=(\mathcal T;
    \ W_Q:=Q^{-1}W,
    \ C_Q:=C,
    \ F_Q:=\left\{p\in F\left\lvert (\ref{HPB of QFib}) \text{ a homotopy pullback} \right.\right\})
$$
\begin{equation}\label{HPB of QFib}\xymatrix@R=0.5cm{
    E\ar@{->>}[d]_p\ar[r]^(0.25){i_E}&E\sqcup_{\overline{\mathscr C} E}QE\ar[d]^{(p,Qp)}\\
    B\ar[r]_(0.25){i_B}&B\sqcup_{\overline{\mathscr C} B}QB
}\end{equation}
The resulting homotopy category is the reflective subcategory
$$\mathcal Ho \mathcal T_Q:=\{X\in \mathcal Ho \mathcal T\mid (i_X:X\rightarrow X\sqcup_{\overline{\mathscr C} X}QX)\in W\}
$$
of $Q$-fibrant objects.
\end{theorem}

The above definition can be dualized and the resulting idempotent quasicomonads induce right Bousfield localizations and associated coreflective homotopy subcategories.

\subsection{Idempotent quasiadjunctions}

A quasiadjunction $(S\dashv_{\ \mathscr C,\mathscr F}\Lambda):\mathcal T\rightleftharpoons \mathcal A$ induces the following natural span on $\mathcal T$ and natural cospan on $\mathcal A$:
$$\xymatrix@C=1.4cm{
    Id_{\mathcal T}
    &\mathscr C\mathfrak C
    \ar@{=>}[l]_{\text{cof}\eta'_{\mathfrak C}}^\sim
    \ar@{=>}[r]^(0.45){(\Lambda\text{fib}_S\eta)_{\mathfrak C }}
    &\Lambda\mathfrak F S\mathfrak C}
    \hspace{0.35cm}
    \xymatrix@C=1.4cm{S\mathfrak C\Lambda\mathfrak F
    \ar@{=>}[r]^{(\epsilon S\text{cof}_\Lambda)_{\mathfrak F }}
    &\mathscr F\mathfrak F
    &Id_{\mathcal A}\ar@{=>}[l]_{\epsilon'_{\mathfrak F}\text{fib}}^\sim
}$$

\begin{definition}
    An \textit{idempotent quasiadjunction} is a quasiadjunction such that the induced span and cospan are respectively an idempotent quasimonad and an idempotent quasicomonad.
\end{definition}

\begin{theorem}[{\cite[Theorem 2.3.8]{Vi20}}]
An idempotent quasiadjunction $(S\dashv_{\ \mathscr C,\mathscr F}\Lambda):\mathcal T\rightleftharpoons \mathcal A$ induces an equivalence between the associated (co)reflective homotopy subcategories.
$$
    \xymatrix{\mathcal Ho \mathcal T\ar@<0.15cm>[r]^(0.4){\mathbb L Id}\ar@{{}{ }{}}[r]|(0.4)\perp&\mathcal Ho \mathcal T_{\Lambda \mathfrak F S\mathfrak C}\ar@<0.15cm>[r]^{\mathbb L S}\ar@{{}{ }{}}[r]|\perp\ar@{^{(}->}@<0.15cm>[l]^(0.6){\mathbb R Id}&\mathcal Ho \mathcal A_{S\mathfrak C\Lambda\mathfrak F}\ar@{^{(}->}@<0.15cm>[r]^(0.6){\mathbb L Id}\ar@{{}{ }{}}[r]|(0.6)\perp\ar@<0.15cm>[l]^{\mathbb R\Lambda}&\mathcal Ho \mathcal A\ar@<0.15cm>[l]^(0.4){\mathbb R Id}}
$$
\end{theorem}

\section{\texorpdfstring{$E^\shortrightarrow_\infty$}{}-algebras}

\subsection{Relative sets and filtered rooted relative trees}

Relative operads are abstract operations with entries indexed by relative sets. We now give the basic definitions and constructions on these colored sets. We will also require filtered rooted relative trees in the construction of the bar resolutions and delooping spectra, and we provide here the relevant definitions and constructions.

Let $\texttt{Set}_{\{d,c\}}$ be the \textit{category of relative sets} composed of sets equipped with a coloring on the colors $\{d,c\}$, ie the class of objects
$$
    \{(A,\mathfrak c)\in\Sigma_{\texttt{Set}}\texttt{Set}(\{A\}\sqcup A,\{d,c\})\mid \mathfrak cA=d\implies \forall a\in A:\mathfrak c a=d\},
$$
with $(A,\mathfrak{c})$ usually being denoted simply as $A$ or explicitly as a set of elements in brackets with coloring given by subscripts, eg $\{1_d,2_d,3_c,4_d,5_c\}_c$. The morphisms sets are 
$$
    \texttt{Set}_{\{d,c\}}(A,A'):=\begin{cases}
        \{\sigma\in\texttt{Set}(A,A')\mid \mathfrak c a=c \implies \mathfrak c' \sigma a=c\},&\!\!\!\!\mathfrak c A=d\text{ or }\mathfrak c'A'=c\\
        \emptyset,&\!\!\!\!\!\mathfrak c A=c\text{ and }\mathfrak c'A'=d
    \end{cases}
$$

For $\star\in \{d,c\}$ we denote by $\texttt{Set}_\star\subset\texttt{Set}_{\{d,c\}}$ the full subcategory of relative sets $A$ such that $\mathfrak{c}A=\star$.

Given $((A,\mathfrak c),\langle (B^a,\mathfrak c^a)\rangle)\in\Sigma_{\texttt{Set}_{\{d,c\}}}\Pi_A\texttt{Set}_{\mathfrak{c}a}$ we have the dependent sum
$$
    (\Sigma_AB^a,\Sigma_A\mathfrak{c}^a)\in \texttt{Set}_{\{d,c\}},\ \Sigma_A\mathfrak{c}^a(\Sigma_AB^a)=\mathfrak{c}A,\ \Sigma_A\mathfrak{c}^a(a,b)=\mathfrak{c}^ab.
$$

For $\sigma\in\texttt{Set}_{\{d,c\}}(A,A')$ let
\begin{gather*}
\sigma(B^a)\in\texttt{Set}_{\{d,c\}}(\Sigma_AB^a,\Sigma_{A'}B^{\sigma^{-1}a'});
\ \sigma(B^a)(a',b):=(\sigma a',b)
\end{gather*}
and for $\langle\tau^a\rangle\in\Pi_A\texttt{Set}_{\mathfrak{c}(a)}(B^a,B^{\prime a})$ let
\begin{gather*}
    \Sigma_A\tau^a\in\texttt{Set}_{\{d,c\}}(\Sigma_AB^a,\Sigma_AB^{\prime a});
    \ \Sigma_A\tau^a(a',b):=(a',\tau^{a'}b).
\end{gather*}

We also have the dependent product
\begin{gather*}
 \Pi_AB^a\in\texttt{Set}_{\{d,c\}},
 \ \Pi_A\mathfrak{c}^a \Pi_AB^a =\mathfrak{c} A ,
 \ \Pi_A\mathfrak{c}^a\langle b^a\rangle=\begin{cases}
	d,&\forall a\in A:\mathfrak{c}^a b^a =d;\\
    c,&\exists a\in A:\mathfrak{c}^a b^a =c.
\end{cases}
\end{gather*}

For $\sigma\in\texttt{Set}_{\{d,c\}}(A,A')$ let
\begin{gather*}
    \sigma\langle B^a\rangle\in\texttt{Set}_{\{d,c\}}(\Pi_AB^a,\Pi_{A'}B^{\sigma^{-1}a'});\ \sigma\langle B^a\rangle\langle b^a\rangle:=\langle b^{\sigma^{-1}a'}\rangle,
\end{gather*}
and for $\langle\tau^a\rangle\in\Pi_A\texttt{Set}_{\mathfrak{c} a}(B^a,B^{\prime a})$ let
\begin{gather*}
    \Pi_A\tau^a\in\texttt{Set}_{\{d,c\}}(\Pi_AB^a,\Pi_AB^{\prime a});\ \Pi_A\tau^a\langle b^a\rangle:=\langle\tau^ab^a\rangle.
\end{gather*}
For every $\langle b^a\rangle\in\Pi_AB^a$ we can form a new relative set $A_{\langle b^a\rangle}$ composed of the pairs $(a,b^a)$ with coloring $\mathfrak{c}_{\langle b^a\rangle}A_{\langle b^a\rangle}=\Pi_A\mathfrak{c}^a\langle b^a\rangle$ and $\mathfrak{c}_{\langle b^a\rangle}(a,b^a)=\mathfrak{c}^ab^a$. This relative set is naturally equipped with $\pi_{\langle b^a\rangle}\in\texttt{Set}_{\{d,c\}}(A_{\langle b^a\rangle},A)$ with $\pi_{\langle b^a\rangle}(a,b^a)=a$. Let
\begin{gather*}
\nu\in\texttt{Set}_{\{d,c\}}(\Pi_A\Sigma_{B^a}C^{a,b},\Sigma_{\Pi_AB^a}\Pi_{A_{(b^a)}}C^{a,b^a});
\hspace{0.5cm}\nu\langle(b^a,c^a)\rangle:=(\langle b^a\rangle,\langle c^a\rangle).
\end{gather*}
This is a key element in distributivity properties.

Let $\mathbb S_{\{d,c\}}^{\text{inj}}\subset \texttt{Set}_{\{d,c\}}$ be the subcategory of $\texttt{Set}_{\{d,c\}}$ composed of the finite relative sets and the injective functions that preserve coloring, i.e.
$$
    \mathbb S_{\{d,c\}}^{\text{inj}}(A,A') = \begin{cases}
        \{\sigma\in\texttt{Set}_{\{d,c\}}(A,A')\mid \sigma \text{ is injective},\ \mathfrak c' \sigma a=\mathfrak c a\},& \mathfrak c A=\mathfrak c' A'\\
        \emptyset,& \mathfrak c A\neq\mathfrak c' A'
    \end{cases}
$$
Let $\mathbb S_{\{d,c\}}\subset \mathbb S_{\{d,c\}}^{\text{inj}}$ be the subcategory with the same objects and bijections that preserve coloring as morphisms. For $\star\in \{d,c\}$ we denote by $\mathbb S^{\text{inj}}_\star$ and $\mathbb S_\star$ the full subcategories of $\mathbb S^{\text{inj}}_{\{d,c\}}$ and $\mathbb S_{\{d,c\}}$ respectively composed of relative sets $A$ such that $\mathfrak{c}A=\star$. Define also the subcategory $\mathbb S_{\langle d<c\rangle}\subset\texttt{Set}_{\{d,c\}}$ with objects the finite relative sets and with morphisms the bijections (that don't necessarily preserve coloring). Note that $\mathbb S_{\{d,c\}}$ is a subcategory of $\mathbb S_{\langle d<c\rangle}$.

Many spaces of interest are built via the two sided bar construction for monads induced by operads, which can be described using filtered rooted relative trees.

\begin{definition} The simplicial category $\mathbb T_{\{d,c\}}\in\texttt{Cat}^{\Delta^{\text{op}}}$ of filtered rooted relative trees has as objects quintuples 
$$T=(\langle V^i\rangle,\langle E^i\rangle,\langle s_i\rangle,\langle t_i\rangle,\mathfrak c)\in \mathbb T_{\{d,c\}}\langle m\rangle$$
composed of
    \begin{itemize}
        \item[-] A sequence of nonempty finite sets $\langle V^{i}\rangle\in \mathbb S^{\langle m-1\rangle}$. We also set $V^{-1}:=\{v^r\}$, and call $v^r$ the \textit{root vertex} of $T$. We set $V:=\sqcup_{\langle m-1\rangle}V^{i}$ and $V_\ast:=V^{-1}\sqcup V$. For $v\in V$ we will denote by $\lvert v\rvert\in\langle m-1\rangle$ the element such that $v\in V^{\lvert v\rvert}$.
        
        \item[-] A sequence of finite sets $\langle E^i\rangle\in  \mathbb S^{\langle  m\rangle}$. We also set $E^{-1}:=\{e^r\}$, and call $e^r$ the \textit{root edge} of $T$. We set $E:=E^{-1}\sqcup(\sqcup_{\langle m\rangle}E^i)$. The edges in $E^m$ are called the \textit{leaves} of $T$.
        
        \item[-] A sequence of bijections $\langle s_{i}\rangle\in\Pi_{\langle  m-1\rangle}\mathbb S(E^{i},V^{i})$, called the \textit{start} of the edges. We sometimes omit the subscript and write simply $se:=s_ie$. Note that the leaves don't have a source.
        
        \item[-] A sequence of functions $\langle t_i\rangle\in\Pi_{\langle m\rangle}\mathbb S(E^i,V^{i-1})$, the \textit{target} of the edges.  We sometimes omit the subscript and write simply $te:=t_ie$. Note that the root edge doesn't have a target.
        
        \item[-] A function $\mathfrak c\in \texttt{Set}(E,\{d,c\})$, the \textit{coloring} of the edges, such that if $te'=se$ and $\mathfrak c e=d$ then $\mathfrak ce'=d$.
    \end{itemize}
    
    We sometimes just write $T=(\langle V^i\rangle,\langle E^i\rangle)$ and leave the mappings implicit.
    
    Morphisms $\sigma\in \mathbb T_{\{d,c\}}\langle m\rangle(T,T')$ are pairs of sequences of bijections
    $$
    (\langle\sigma^i_V\rangle,\langle\sigma^i_E\rangle)\in\Pi_{\langle m-1\rangle}\mathbb S^{\text{inj}}(V^i,V^{\prime,i})\times \Pi_{\langle m\rangle}\mathbb S^{\text{inj}}(E^i,E^{\prime,i})
    $$
    that commute with the structural functions.
    
    The simplicial structural functors are defined on objects as:
    \begin{align*}
        T\cdot \partial_i&:=\left(\langle V^{\partial_ij}\rangle,\langle E^{\partial_ij}\rangle, \langle s_{\partial_ij}\rangle,\left\langle \begin{cases}
           t_{\partial_ij},&j\neq i\\
           t_is_i^{-1}t_{i+1},&j=i.
        \end{cases}\right\rangle,\mathfrak c\cdot \partial_i\right)\\
        T\cdot \delta_i&:=\left(\langle V^{\delta_ij}\rangle,\langle E^{\delta_ij}\rangle, \langle s_{\delta_ij}\rangle,\left\langle \begin{cases}
           t_{\delta_ij},&j\neq i+1\\
           s_i,&j=i+1.
        \end{cases}\right\rangle,\mathfrak c\cdot \delta_i\right)
    \end{align*}
    with the coloring maps induced naturally from the ones in $T$.
    
    Define also $\mathbb T_{\{d,c\}}^0\in\texttt{Cat}^{\Delta^{\text{op}}}$ as the full simplicial subcategory of relative trees such that $|V^1|=1$. Define also the simplicial full subcategories $\mathbb T_\star\subset \mathbb T_{\{d,c\}}$ for $\star\in\{d,c\}$ of the trees such that $\mathfrak c e^r=\star$. We similarly define the simplicial full subcategories $\mathbb T_\star^0\subset \mathbb T_{\{d,c\}}^0$.
    
    Note that any $T\in\mathbb T_{\{d,c\}}\langle m\rangle$ has a natural partial order structure on the union of the set of vertices and edges induced by the start and target maps such that $e^r$ is the unique minimal element. For each $e\in E^1$ let $T_{\geq e}\in \mathbb T_{\{d,c\}}^0\langle m\rangle$ be the sub-tree composed of the root vertex, root edge and the vertices and edges greater than $e$.
    
    For all $T\in\mathbb T_{\{d,c\}}\langle m\rangle$ and $v\in V$ define the relative set
    $$\text{in }v:=\{e\in E\mid t e=v\}_{\mathfrak c s^{-1}v}\in\mathbb S_{\{d,c\}}.
    $$
\end{definition}

Note that $\mathbb S_{\{d,c\}}$ is isomorphic to $\mathbb T_{\{d,c\}}\langle 0\rangle$.

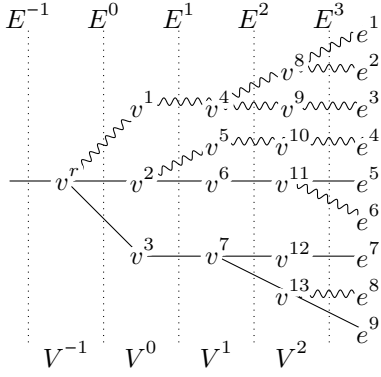
\begin{figure}
    \centering
    \begin{tikzpicture}
        \draw [decorate,decoration={snake,segment length=1.5mm,amplitude=0.5mm}] (0,0) -- (1,1) -- (2,1) -- (3,1.5)--(4,2);
        \draw[decorate,decoration={snake,segment length=1.5mm,amplitude=0.5mm}] (3,1.5)--(4,1.5);
        \draw[decorate,decoration={snake,segment length=1.5mm,amplitude=0.5mm}] (2,1) -- (3,1)--(4,1);
        \draw [decorate,decoration={snake,segment length=1.5mm,amplitude=0.5mm}]
        (1,0)--(2,0.5)--(3,0.5)--(4,0.5);
        \draw (-0.75,0)--(1,0)--(2,0)--(3,0)--(4,0);
        \draw[decorate,decoration={snake,segment length=1.5mm,amplitude=0.5mm}] (3,0)--(4,-0.5);
        \draw (0,0) -- (1,-1) -- (2,-1) -- (3,-1)--(4,-1);
        \draw[decorate,decoration={snake,segment length=1.5mm,amplitude=0.5mm}] (3,-1.5)--(4,-1.5);
        \draw (2,-1) -- (3,-1.5)--(4,-2);
        
        \draw (0.025,0.025) node {\contour{white}{$v^r$}};
        
        \draw (1.02,1.02) node {\contour{white}{$v^{1}$}};
        \draw (1.02,0.02) node {\contour{white}{$v^{2}$}};
        \draw (1.02,-0.93) node {\contour{white}{$v^{3}$}};
        
        \draw (2.02,1.02) node {\contour{white}{$v^4$}};
        \draw (2.02,0.52) node {\contour{white}{$v^5$}};
        \draw (2.02,0.02) node {\contour{white}{$v^6$}};
        \draw (2.02,-0.93) node {\contour{white}{$v^{7}$}};
        
        \draw (3.02,1.52) node {\contour{white}{$v^8$}};
        \draw (3.02,1.02) node {\contour{white}{$v^9$}};
        \draw (3.02,0.52) node {\contour{white}{$v^{10}$}};
        \draw (3.02,0.02) node {\contour{white}{$v^{11}$}};
        \draw (3.02,-0.98) node {\contour{white}{$v^{12}$}};
        \draw (3.02,-1.48) node {\contour{white}{$v^{13}$}};

        \draw (4.02,2.02) node {\contour{white}{$e^1$}};
        \draw (4.02,1.52) node {\contour{white}{$e^2$}};
        \draw (4.02,1.02) node {\contour{white}{$e^3$}};
        \draw (4.02,0.52) node {\contour{white}{$e^4$}};
        \draw (4.02,0.02) node {\contour{white}{$e^5$}};
        \draw (4.02,-0.48) node {\contour{white}{$e^6$}};
        \draw (4.02,-0.98) node {\contour{white}{$e^7$}};
        \draw (4.02,-1.48) node {\contour{white}{$e^8$}};
        \draw (4.02,-1.98) node {\contour{white}{$e^9$}};
        
        \draw[dotted] (-0.5,2)--(-0.5,-2.2);
        \draw (-0.5,2.2) node {$E^{-1}$};
        \draw[dotted] (0.5,2)--(0.5,-2.2);
        \draw (0.5,2.2) node {$E^0$};
        \draw[dotted] (1.5,2)--(1.5,-2.2);
        \draw (1.5,2.2) node {$E^1$};
        \draw[dotted] (2.5,2)--(2.5,-2.2);
        \draw (2.5,2.2) node {$E^2$};
        \draw[dotted] (3.5,2)--(3.5,-2.2);
        \draw (3.5,2.2) node {$E^3$};
        
        \draw (0,-2.3) node {$V^{-1}$};
        \draw (1,-2.3) node {$V^0$};
        \draw (2,-2.3) node {$V^1$};
        \draw (3,-2.3) node {$V^2$};
        
    \end{tikzpicture}
    \caption{A filtered rooted relative tree in $\mathbb T_c\langle 3\rangle$ with wiggled edges representing ``domain'' edges and straight edges ``codomain'' edges. The leaves are the only edges that are labeled.}
\end{figure}

We have natural dependent sums and dependent products of filtered rooted trees of a fixed height $\langle T^a\rangle\in\Pi_A\mathbb T_{\{d,c\}}\langle m\rangle$ defined as
\begin{gather*}
    \Sigma_A T^a:=(\langle \Sigma_A V^{a,i}\rangle,
    \langle \Sigma_A E^{a,i}\rangle,
    \langle \Sigma_A s_i^a\rangle,
    \langle \Sigma_A t_i^a\rangle,
    \Sigma_A\mathfrak c^a),\\
    \Pi_A T^a:=(\langle \Pi_A V^{a,i}\rangle,
    \langle \Pi_A E^{a,i}\rangle,
    \langle \Pi_A s_i^a\rangle,
    \langle \Pi_A t_i^a\rangle,
    \Pi_A\mathfrak c^a).
\end{gather*}

We also have for all 
$$
    T=(\langle V^i\rangle, \langle E^i\rangle)\in\mathbb T_{\{d,c\}}\langle m\rangle
    \text{ and }
    \langle S^e\rangle=\langle(\langle W^{e,i}\rangle, \langle F^{e,i}\rangle)\rangle\in\Pi_{E^m}\mathbb T_{\mathfrak c e}\langle n\rangle
$$
the grafting 
\begin{gather*}
    T\circ\langle S^e\rangle
    :=\left(\left\langle \begin{cases}
    V^i,& i< m\\
    \Sigma_{E^m}W^{e,i-m-1},&i\geq m
    \end{cases}
    \right\rangle,\  
    \left\langle \begin{cases}
    E^i,& i\leq m\\
    \Sigma_{E^m}F^{e,i-m-1},&i>m
    \end{cases}
    \right\rangle\right)
\end{gather*}
in $\mathbb T_{\{d,c\}}\langle m+n+1\rangle$, with the obvious start, target and coloring maps.

\subsection{Relative operads}

We now give a brief review of relative operads, a kind of colored operad introduced by Voronov in \cite{Vo99}.

\begin{definition}\label{RelOperad}
The category of $\mathbb S_{\{d,c\}}$-spaces is the contravariant functor category $\texttt{Top}^{\mathbb S_{\{d,c\}}^{\texttt{op}}}$. A topological relative operad is an $\mathbb S_{\{d,c\}}$-space $\mathcal P\in \texttt{Top}^{\mathbb S_{\{d,c\}}^{\texttt{op}}}$ equipped with elements $\textbf{id}_\star\in\mathcal P\{1_\star\}_\star$ for $\star\in\{d,c\}$ and structural maps
    $$
        \langle \circ_{A,\langle B^a\rangle}\rangle\in\Pi_{\Sigma_{\mathbb S_{\{d,c\}}}\Pi_A\mathbb S_{\mathfrak c a}}\texttt{Top}(\mathcal PA\times\Pi_A\mathcal PB^a,\mathcal P\Sigma_AB^a)
    $$
such that $\mathcal P\emptyset_\star=\ast$ for $\star\in\{d,c\}$ and, using the notation
    $$
        \alpha\left\langle \beta^a\right\rangle:=\circ_{A,\langle B^a\rangle}\left(\alpha,\left\langle \beta^a\right\rangle\right),
    $$
the following equations are satisfied:
    \begin{align*}
        \alpha\langle \beta^a \langle \gamma^{a,b}\rangle\rangle&=\alpha \langle \beta^a\rangle\langle \gamma^{a,b}\rangle;\\
        \textbf{id}_{\mathfrak{c}A} \alpha=&\ \alpha=\alpha\left\langle \textbf{id}_{\mathfrak{c}a}\right\rangle;\\
        \alpha\cdot\sigma\langle \beta^a\rangle&=\alpha\langle \beta^{\sigma^{-1}a'}\rangle\cdot\sigma(B^a);\\
        \alpha\langle \beta^a\cdot\tau^a\rangle&=\alpha\langle \beta^a\rangle\cdot\Sigma_A\tau^a.
    \end{align*}
    
    Operad morphisms are natural transformations that preserve the unit and compositions, and we denote the category of topological relative operad as $\texttt{Op}_{\{d,c\}}[\texttt{Top}]$.
\end{definition}

For $X=((X_d,e_d),(X_c,e_c))\in\texttt{Top}^2_\ast$ define 
$$
\Pi_-X:\mathbb S^{\text{inj}}_{\{d,c\}}\rightarrow\texttt{Top};\ \Pi_A X:=\Pi_A X_{\mathfrak{c}a},\ \sigma\cdot \langle x^a\rangle:=\left\langle\begin{cases}
        e_{\mathfrak{c}a'},&a'\not\in\text{Im }\sigma;\\
        x^{\sigma^{-1}a'},&a'\in\text{Im }\sigma.
    \end{cases}\right\rangle.
$$

The underlying functor of a unital relative operad $\mathcal P$ can be extended to a functor on $\mathbb S_{\{d,c\}}^{\text{inj},\text{op}}$. For $\sigma\in\mathbb S_{\{d,c\}}^{\text{inj}}(A,A')$ the right action $\cdot\sigma\in\texttt{Top}(\mathcal PA',\mathcal PA)$ is defined as 
$$
\alpha\cdot\sigma:=\alpha\left\langle\begin{cases}
\ast_{\mathfrak c a'},&a'\not\in\text{Im }\sigma;\\
\textbf{id}_{\mathfrak{c}a'},&a'\in\text{Im }\sigma.
\end{cases}\right\rangle.
$$
These morphisms are the \textit{degenerations} of the relative operad.

A relative operad $\mathcal P$ induces a monad $(P;\eta,\mu)$ on $\texttt{Top}_\ast^2$ with
\begin{gather*}
    PX_\star:=\textstyle  \int^{\mathbb S_\star^{\text{inj}}}\mathcal P A \times \Pi_A X;\\
    \eta_\star x:=[\textbf{id}_\star,x],
    \hspace{1cm}\mu_\star[\alpha,\langle[\beta^a,\langle x^{a,b}\rangle]\rangle]:=[\alpha \langle \beta^a\rangle,\langle x^{a,b}\rangle].
\end{gather*}

\begin{definition}\label{P-alg}
    Let $\mathcal P$ be a relative operad. A \textit{$\mathcal P$-space} is a $P$-algebra, i.e. a pair of pointed spaces $X\in\texttt{Top}^2_\ast$ equipped with structural maps
    $$
        \langle \theta_A\rangle \in \Pi_{\mathbb S_{\{d,c\}}}\texttt{Top}(\mathcal P A \times \Pi_A X,X_{\mathfrak c A}),
    $$
	satisfying, using the notation $\alpha\langle x^a\rangle=\theta_A\langle \alpha,\langle x^a\rangle\rangle$, the following equations:
\begin{align*}
\alpha\left\langle \beta^a\left\langle x^{a,b}\right\rangle\right\rangle&=\alpha\circ\langle \beta^a\rangle\left\langle x^{a,b}\right\rangle;\\
\textbf{id}_\star x&=x;\\
\alpha\cdot\sigma\langle x^a\rangle&=\alpha( \sigma\cdot\langle x^a\rangle).
\end{align*}
The category of $\mathcal P$-spaces is denoted $\mathcal P[\texttt{Top}]$.
\end{definition}

The following are the relative operads relevant to the main result.

The terminal relative operad is $Com^\shortrightarrow$ with underlying $\mathbb S_{\{d,c\}}$-space given by $Com^\shortrightarrow(A):=\ast$. The $\mathbb S_{\{d,c\}}$ right actions, units and compositions are the unique terminal maps. The $Com^\shortrightarrow$-spaces are pairs $(M_d,M_c)$ of topological commutative monoids equipped with a continuous homomorphism $\iota:M_d\rightarrow M_c$ induced by the unique element in $Com^\shortrightarrow\{1_c\}_o$.

For $U\in\mathscr A$ let the relative operad of $U$-embeddings $\text{Emb}_U^{\shortrightarrow}$ is
\begin{gather*}
    \text{Emb}_U^{\shortrightarrow}A:=\{\alpha=\langle\alpha_a\rangle\in U^{\sqcup_A U}\mid \langle\alpha_a\rangle\text{ is an embedding}\};\\
    \langle \alpha_{a'}\rangle\cdot\sigma:=\langle \alpha_{\sigma a}\rangle,
    \hspace{0.5cm}  \textbf{id}_\star:=id_U,
    \hspace{0.5cm} \alpha\langle \beta^a\rangle:=\langle \alpha_a\beta^a_b\rangle
\end{gather*}
and degenerations deleting embeddings.

For $U\in\mathscr A$ the loop space map functors image has natural $\text{Emb}_U^\shortrightarrow$-pairs structure, giving us the functor
\begin{gather}
    \nonumber\Omega^U_2:\texttt{Top}_\ast^\shortrightarrow\rightarrow \text{Emb}_U^\shortrightarrow[\texttt{Top}];
    \hspace{0.3cm}\Omega^U_2(\iota:Y_d\rightarrow Y_c):=(Y_d^{\mathds S^U},Y_c^{\mathds S^U});\\
    \label{SumRelLopSpc}\alpha\langle \gamma^a\rangle:=\left(\vec u\mapsto\begin{cases}
        \gamma^a\alpha_a^{-1}\vec u,&\mathfrak{c}a=\mathfrak{c}A\\
        \iota\gamma^a\alpha_a^{-1}\vec u,&\mathfrak{c}a\neq\mathfrak{c}A
    \end{cases}\right)
\end{gather}
For $(U,V)\in\Sigma_{\mathscr A}\mathscr A_U$ we have natural inclusion of relative operads
\begin{equation*}
    i^U_V:\text{Emb}_U^\shortrightarrow\Rightarrow \text{Emb}_V^\shortrightarrow;
    \hspace{0.5cm}i^U_V \alpha := \langle \vec v\mapsto\vec v_{V-U}+\alpha_a\vec v_U\rangle
\end{equation*}
and we define $\text{Emb}^\shortrightarrow_\infty:=\text{colim}_{\mathscr A}\text{Emb}^\shortrightarrow_U$.

The embeddings operad contains embeddings of  configuration spaces, and these embeddings are relevant to the definition of $E^\shortrightarrow$-operads we give here. For each $U\in\mathscr A$ define the configurations $\mathbb S_{\{d,c\}}$-space
$$
    \text{Conf}_U^\shortrightarrow:\mathbb S_{\{d,c\}}^{\text{op}}\rightarrow\texttt{Top};\ \text{Conf}_U^\shortrightarrow A:=\{\vec x=\langle \vec x_a\rangle\in U^A\mid a\neq a'\implies \vec x_{a}\neq \vec x_{a'}\}.
$$

Note that $\text{Conf}_U^\shortrightarrow$ is $m$-cofibrant, since it is $h$-equivalent to the underlying space of the Fulton-MacPherson operads which are $q$-cofibrant \cite{Sa01,Ho12}. We can define the $\mathbb S_{\{d,c\}}$-space maps
$$
\chi_U:\text{Conf}_U^\shortrightarrow\Rightarrow \text{Emb}_U^\shortrightarrow;\ \chi_U\vec x:=\left\langle\vec u\mapsto\vec x_a+\frac{\min_{a'\neq a''}\lVert \vec x_{a'}-\vec x_{a''}\rVert\vec u}{\min_{a'\neq a''}\lVert \vec x_{a'}-\vec x_{a''}\rVert+2\lVert\vec u\rVert}\right\rangle.
$$

\begin{definition}
    An $E_\infty^\shortrightarrow$-operad is an operad 
    $$
        \mathcal E^\shortrightarrow\in\texttt{Op}_{\{d,c\}}[\texttt{Top}]
    $$
    equipped with a relative operad map 
    $$
        \Psi\in\texttt{Op}_{\{d,c\}}[\texttt{Top}](\mathcal E^\shortrightarrow,\text{Emb}^\shortrightarrow_\infty)
    $$
    and, for the induced $\mathscr A$-filtration $\mathcal E^\shortrightarrow_U:=\Psi^{-1}\text{Emb}_U$, a $\mathbb S_{\{d,c\}}$-space homotopy equivalence
    $$
        \Phi_U\texttt{Top}^{\mathbb S_{\{d,c\}}^{\text{op}}}(\text{Conf}_U,\mathcal E^\shortrightarrow_U)
    $$
    for each $U\in\mathscr A$ such that $\Psi\restriction_U\Phi_U=\chi_U$.
\end{definition}

By this definition the $\mathcal E^\shortrightarrow_U$ are $m$-cofibrant as $\mathbb S_{\{d,c\}}$-spaces and $\mathcal E^\shortrightarrow$ is contractible and free. One of the main examples of $E_\infty^\shortrightarrow$-operads we will consider is the Steiner relative operad, composed of paths of embeddings \cite{St79}.

For all $U\in\mathscr A$ define the relative operad $\mathscr H^{\shortrightarrow}_U$ as
\begin{gather*}
    \mathscr H^{\shortrightarrow}_UA:= \left\{\alpha=\langle\alpha_a\rangle\in U^{\sqcup_A I\times U}\left\lvert\begin{matrix*}[l]
    \forall a\in A, t\in I:(\vec u\mapsto\alpha_a(t,\vec u))\in \text{Emb}_U^\shortrightarrow\{a\};\\
    \forall a\in A,\forall t\in I,\forall \vec u,\vec v\in U:\\
    \lVert \alpha_a(t,\vec u)- \alpha_a(t,\vec v)\rVert\leq \lVert\vec u- \vec v\rVert;\\
    \forall a\in A,\vec u\in U:\alpha_a(1,\vec u)=\vec u;\\
    \langle \vec u\mapsto\alpha_a(0,\vec u)\rangle\in \text{Emb}_U^\shortrightarrow A.
    \end{matrix*}\right.\right\}\\
    \langle \alpha_{a'}\rangle\cdot\sigma:=\langle\alpha_{\sigma a}\rangle,\ 
    \textbf{id}_\star:=((t,\vec u)\mapsto \vec u),\ 
    \alpha \langle \beta^a\rangle:=\langle ( t,\vec u)\mapsto\alpha_a(t,\beta^a_b(t,\vec u)) \rangle
\end{gather*}
and degenerations deleting paths of embeddings.

We have natural inclusions $\iota^U_V:\mathscr H^\shortrightarrow_U\Rightarrow\mathscr H^\shortrightarrow_V$ for all $(U,V)\in\Sigma_{\mathscr A}\mathscr A_U$ with
$$
    \iota^U_V \alpha :=\langle (t,\vec v)\mapsto\vec v_{V-U}+\alpha_a(t,\vec v_U)\rangle
$$
and we define $\mathscr H^\shortrightarrow_\infty:=\text{colim}_{\mathscr A}\mathscr H^\shortrightarrow_U$.

The $E_\infty^\shortrightarrow$-structural transformations are
\begin{gather*}
    \Psi:\mathscr H^{\shortrightarrow}\Rightarrow\text{Emb}^\shortrightarrow;\ \Psi_U\alpha:=\langle \vec u\mapsto\alpha_a(0,\vec u)\rangle;\\
    \Phi_U:  \text{Conf}_U^\shortrightarrow \Rightarrow \mathscr H_U^\shortrightarrow;
    \ \Phi_U\vec x:=\left\langle (t,\vec u)\mapsto(1-t)\left((\chi_U\vec x)_a\vec u\right)+ t\vec u\right\rangle.
\end{gather*}
The homotopy inverses of the $\Phi_U$ are
$$
    \bar\Phi_U:\mathscr H^\shortrightarrow_U\Rightarrow\text{Conf}^\shortrightarrow_U;\ \bar\Phi_U\alpha=\langle \alpha_a(0,\vec 0) \rangle.
$$
See \cite{St79} for the construction of the homotopies.

\subsection{Bar resolution}

For the construction of the quasiadjunctions in our main theorems we will require the \textit{bar resolution} of $\mathcal E^\shortrightarrow$-pairs. Recall from \cite[Construction 9.6]{MayGeomItLoopSp} that for a monad $(C,\eta,\mu)$ in the category $\mathcal T$, a $C$-functor $(F,\lambda)$ in the category $\mathcal A$ and a $C$-algebra $(X,\xi)$ the two sided bar construction $B_-(F,C,X)\in \mathcal A^\Delta$ with
\begin{gather*}
    B_{\langle m\rangle}(F,C,X):=FC^m X;
    \delta_i:=FC^i\eta_{C^{m-i}},
    \partial_i:=\begin{cases}
        \lambda_{C^m},&i=0;\\
        FC^{i-1}\mu_{C^{m-i+1}},&0<i<m;\\
        FC^{m-1}\xi,&i=m.
    \end{cases}
\end{gather*}

In particular for a relative operad $\mathcal P$ and $C= P=F$ we have a natural isomorphism
\begin{gather*}
    \textstyle B_{\langle m\rangle}(P,P,X)_\star\cong
    \int^{\mathbb T_\star\langle m\rangle}\Pi_{V_\ast}\mathcal P\text{in }v\times \Pi_{E^m}X_{\mathfrak ce},\\
    [\alpha^r,\langle \alpha^v\rangle,\langle x^e\rangle]_T\cdot\partial_i:=\\
    \begin{cases}
        [\alpha^r\langle\alpha^{se'}\rangle,\langle \alpha^v\rangle,\langle x^e\rangle]_{T\cdot\partial_0},&i=0\\
        \left[\alpha^r,\left\langle \begin{cases}
            \alpha^v\langle \alpha^{se'}\rangle,
            &\lvert v\rvert=i-1\\
            \alpha^v,
            &\lvert v\rvert\neq i-1
        \end{cases} \right\rangle,\langle x^e\rangle\right]_{T\cdot\partial_i},&0<i<m\\
        [\alpha^r,\langle \alpha^v\rangle,\langle\alpha^{s^{-1}e}\langle x^{e'}\rangle\rangle]_{T\cdot\partial_m},&i=m
    \end{cases},\\
    [\alpha^r,\langle\alpha^v\rangle,\langle x^e\rangle]_T\cdot\delta_i:=
    \left[\alpha^r,
    \left\langle \begin{cases}
        \textbf{id}_{\mathfrak c s^{-1}v}, &\lvert v\rvert=i\\
        \alpha^v, &\lvert v\rvert\neq i
    \end{cases}\right\rangle,
    \langle x^e\rangle
    \right]_{T\cdot\delta_i}.
\end{gather*}

The $\mathcal E^\shortrightarrow$-pair structural maps in each dimension are
\begin{equation}\label{SimplAlphaStrMap}
    \alpha\langle[\beta^{a,r},\langle\beta^{a,v}\rangle,\langle x^{a,e}\rangle]_{T^a}\rangle:=[\alpha\langle\beta^{a,r}\rangle,\langle\beta^{a,v}\rangle,\langle x^{a,e}\rangle]_{\Sigma_A T^a}
\end{equation}

The bar resolution of $\mathcal E^\shortrightarrow$-pairs is then the geometric realization of this simplicial $\mathcal E^\shortrightarrow$-pair functor
\begin{gather*}
    \textstyle\overline B_2:\mathcal E^\shortrightarrow[\texttt{Top}]\rightarrow \mathcal E^\shortrightarrow[\texttt{Top}],
    \hspace{0.5cm}
    \overline B_2X_\star
    :=\lvert B_-(E^\shortrightarrow,E^\shortrightarrow,X)_\star\rvert.
\end{gather*}

By the above isomorphism we can intuitively think of points in $\overline B_2 X$ as equivalence classes of filtered rooted relative trees with vertices decorated with elements of $\mathcal E^\shortrightarrow$, leaves decorated with elements of $X$ and we associate an ordered partition with the filtration of the inner vertices.

\begin{figure}
    \centering
    \begin{tikzpicture}
        
    \draw (0,0) node {$\alpha$};
    
    \draw (0.4,-1) -- (0.2,0) --  (0.4,1);
    
    \draw (0.3,0)--(0.8,0)--(1.6,0.5)--(2.4,0.75);
    \draw (1.6,0.5)--(2.4,0.25);
    \draw (0.8,0)--(1.6,-0.5)--(2.4,-0.75);
    \draw[decorate,decoration={snake,segment length=1.5mm,amplitude=0.5mm}] (1.6,-0.5)--(2.4,-0.25);
        
    \draw (0.8,0) node {\tiny\contour{white}{$\beta^{1,r}$}};
        
    \draw[dashed] (1.6,-1)--(1.6,1);
    \draw (1.6,-1.3) node {\tiny$\frac{1}{3}$};    
        
    \draw (1.6,0.5) node {\tiny\contour{white}{$\beta^{1,1}$}};
    \draw (1.6,-0.5) node {\tiny\contour{white}{$\beta^{1,2}$}};
    
    \draw (2.4,0.75) node {\tiny\contour{white}{$x^{1,1}$}};
    \draw (2.4,0.25) node {\tiny\contour{white}{$x^{1,2}$}};
    \draw (2.4,-0.25) node {\tiny\contour{white}{$x^{1,3}$}};
    \draw (2.4,-0.75) node {\tiny\contour{white}{$x^{1,4}$}};
        
    \draw (2.8,0) node {$,$};
        
    \draw[decorate,decoration={snake,segment length=1.5mm,amplitude=0.5mm}] (3,0)--(3.8,0)--(4.6,0.5)--(5.4,0)--(6.2,0.5);
    \draw[decorate,decoration={snake,segment length=1.5mm,amplitude=0.5mm}] (4.6,0.5)--(5.4,0.5);
    
    \draw[decorate,decoration={snake,segment length=1.5mm,amplitude=0.5mm}] (5.4,0)--(6.2,0);
    \draw[decorate,decoration={snake,segment length=1.5mm,amplitude=0.5mm}] (3.8,0)--(4.6,-0.5)--(5.4,-0.5)--(6.2,-0.5);
        
    \draw (3.8,0) node {\tiny\contour{white}{$\beta^{2,r}$}};
    
    \draw[dashed] (4.6,-1)--(4.6,1);
    \draw (4.6,-1.3) node {\tiny$\frac{1}{6}$};
    
    \draw (4.6,0.5) node {\tiny\contour{white}{$\beta^{2,1}$}};
    \draw (4.6,-0.5) node {\tiny\contour{white}{$\beta^{2,2}$}};
    
    \draw[dashed] (5.4,-1)--(5.4,1);
    \draw (5.4,-1.3) node {\tiny$\frac{2}{3}$};
    
    \draw (5.4,0.5) node {\tiny\contour{white}{$\ast$}};
    \draw (5.4,0) node {\tiny\contour{white}{$\beta^{2,4}$}};
    \draw (5.4,-0.5) node {\tiny\contour{white}{$\beta^{2,5}$}};
    
    \draw (6.2,0.5) node {\tiny\contour{white}{$x^{2,1}$}};
    \draw (6.2,0) node {\tiny\contour{white}{$x^{2,2}$}};
    \draw (6.2,-0.5) node {\tiny\contour{white}{$x^{2,3}$}};
    
    \draw (6.6,-1) -- (6.8,0) --  (6.6,1);

    \draw (6.95,0) node {$=$};
    
    \draw (7.15,0)--(8.3,0)--(9.2,1.25)--(10,1.25)--(10.8,1.5)--(11.6,1.5);
    \draw (10,1.25)--(10.8,1)--(11.6,1);
    \draw[decorate,decoration={snake,segment length=1.5mm,amplitude=0.5mm}] (10,0.25)--(10.8,0.5)--(11.6,0.5);
    \draw (8.3,0)--(9.2,0.25)--(10,0.25)--(10.8,0)--(11.6,0);
    \draw[decorate,decoration={snake,segment length=1.5mm,amplitude=0.5mm}] (8.3,0)--(9.2,-0.5)--(10.8,-0.5);
    \draw[decorate,decoration={snake,segment length=1.5mm,amplitude=0.5mm}] (9.2,-0.5)--(10,-1)--(10.8,-1)--(11.6,-0.5);
    \draw[decorate,decoration={snake,segment length=1.5mm,amplitude=0.5mm}] (10.8,-1)--(11.6,-1);
    
    \draw[decorate,decoration={snake,segment length=1.5mm,amplitude=0.5mm}] (8.3,0)--(9.2,-1.44)--(11.6,-1.44);

    \draw (8.25,0) node {\tiny\contour{white}{$\alpha\langle \beta^{1,r},\beta^{2,r}\rangle$}};
    
    \draw[dashed] (9.2,-1.7)--(9.2,1.7);
    \draw (9.2,-2) node {\tiny$\frac{1}{6}$};
    
    \draw (9.2,1.25) node {\tiny\contour{white}{$\textbf{id}_c$}};
    \draw (9.2,0.25) node {\tiny\contour{white}{$\textbf{id}_c$}};
    \draw (9.2,-0.5) node {\tiny\contour{white}{$\beta^{2,1}$}};
    \draw (9.2,-1.5) node {\tiny\contour{white}{$\beta^{2,2}$}};
    
    \draw[dashed] (10,-1.7)--(10,1.7);
    \draw (10,-2) node {\tiny$\frac{1}{3}$};
    
    \draw (10,1.25) node {\tiny\contour{white}{$\beta^{1,1}$}};
    \draw (10,0.25) node {\tiny\contour{white}{$\beta^{1,2}$}};
    \draw (10,-0.5) node {\tiny\contour{white}{$\textbf{id}_d$}};
    \draw (10,-1) node {\tiny\contour{white}{$\textbf{id}_d$}};
    \draw (10,-1.5) node {\tiny\contour{white}{$\textbf{id}_d$}};
    
    \draw[dashed] (10.8,-1.7)--(10.8,1.7);
    \draw (10.8,-2) node {\tiny$\frac{2}{3}$};
    
    \draw (10.8,1.5) node {\tiny\contour{white}{$\textbf{id}_c$}};
    \draw (10.8,1) node {\tiny\contour{white}{$\textbf{id}_c$}};
    \draw (10.8,0.5) node {\tiny\contour{white}{$\textbf{id}_d$}};
    \draw (10.8,0) node {\tiny\contour{white}{$\textbf{id}_c$}};
    \draw (10.8,-0.5) node {\tiny\contour{white}{$\ast$}};
    \draw (10.8,-1) node {\tiny\contour{white}{$\beta^{2,4}$}};
    \draw (10.8,-1.5) node {\tiny\contour{white}{$\beta^{2,4}$}};

    \draw (11.5,1.5) node {\tiny\contour{white}{$x^{1,1}$}};
    \draw (11.5,1) node {\tiny\contour{white}{$x^{1,2}$}};
    \draw (11.5,0.5) node {\tiny\contour{white}{$x^{1,3}$}};
    \draw (11.5,0) node {\tiny\contour{white}{$x^{1,4}$}};
    \draw (11.5,-0.5) node {\tiny\contour{white}{$x^{2,1}$}};
    \draw (11.5,-1) node {\tiny\contour{white}{$x^{2,2}$}};
    \draw (11.5,-1.5) node {\tiny\contour{white}{$x^{2,3}$}};

    \end{tikzpicture}
    \caption{$E^\shortrightarrow_\infty$-structure of $\overline B_2X$}
    \label{fig:E str of B2X}
\end{figure}
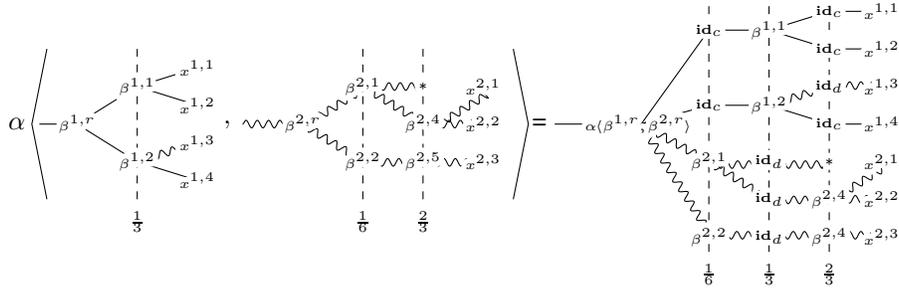

It is not the case in general that the geometric realization of a simplicial $C$-algebra for a topological monad $C$ is a $C$-algebra. This is however the case when the monad is the one induced by an operad. The structural maps are induced by the algorithm described in \cite[Theorem 11.5]{MayGeomItLoopSp}. For a sequence of elements with representatives of distinct dimensions we can systematically determine equivalent representatives of the same dimension, and then apply the formula \ref{SimplAlphaStrMap}, so that the $\mathcal E^\shortrightarrow$-pair structural maps of $\overline B_2X$ are defined by the formula
\begin{multline}\label{EstructBarB}
	\alpha\langle [[\beta^{a,r},\langle \beta^{a,v}\rangle,\langle x^{a,e}\rangle]_{T^a},t^a]\rangle
    :=\\
    \left[\left[\alpha\langle \beta^{a,r}\rangle,\left\langle \begin{cases}
        \textbf{id}_{\mathfrak c s^{-1}v},
        &\exists i>\lvert v\rvert: \delta^ai= \delta^a\lvert v\rvert\\
        \beta^{a,v},
        &\text{otherwise}
    \end{cases}\right\rangle,\langle x^{a,e}\rangle\right]_{\Sigma_AT^a\cdot\delta^a},\lhd_A t^a\right]
\end{multline}
which is illustrated in figure \ref{fig:E str of B2X}.

This functor can be equipped with the natural transformation \begin{equation}\label{eta'}
    \eta':\overline B_2\Rightarrow Id,
    \hspace{1cm} \eta_\star'\left[[\alpha^r,\langle \alpha^v\rangle,\langle x^e\rangle]_T, t\right]
    :=\circ_T\alpha^v\langle x^e\rangle.
\end{equation}
where $\circ_T\alpha^v$ is the composition of all the $\alpha^v$, including $\alpha^r$, induced by the operadic composition and the structure of $T$.

\subsection{Relative operad action}

Operad actions encodes distributive laws between operations defined by operads \cite[Definition VI.1.6]{MaEinftyRingSpcsSpectra}. The following definition is a relative version of this notion.

\begin{definition}\label{RelOpAct}
     A \textit{relative operad pair} is a pair of relative operads $(\mathcal P,\mathcal G)$ equipped with an extension of $\mathcal G$ to $\mathbb S_{\langle d<c\rangle}^{\text{op}}$ and structural maps
    $$
        \textstyle\langle\ltimes_{A,\langle B^a\rangle}\rangle\in\Pi_{\Sigma_{\mathbb S_{\{d,c\}}}\Pi_A\mathbb S_{\mathfrak c a}}\texttt{Top}\left(\mathcal G A \times\Pi_A\mathcal{P}B^a,\mathcal{P}\Pi_AB^a\right)
    $$
the \textit{action} of $\mathcal G$ on $\mathcal P$, such that, using the notations
$$
    f\ltimes\left\langle \alpha^a\right\rangle:=\ltimes\left\langle f,\left\langle \alpha^a\right\rangle\right\rangle,
    \hspace{1cm} f \ltimes\left\langle \beta^{a,b^a}\right\rangle:= f\cdot \pi_{\langle b^a\rangle}  \ltimes\left\langle \beta^{a,b^a}\right\rangle
$$
the following equations are satisfied:
\begin{align}
    f\ltimes\langle g^a\ltimes\langle \alpha^{a,b}\rangle \rangle&
    =f\circ\langle g^a\rangle\ltimes\langle \alpha^{a,b}\rangle;\label{OpActAct}\\
    f\ltimes\langle \alpha^a\circ\langle \beta^{a,b}\rangle\rangle&
    =f\ltimes\langle \alpha^a\rangle\circ\langle f \ltimes\langle \beta^{a,b^a}\rangle\rangle\cdot\nu;\label{OpActDistr}\\
    \textbf{id}_{\mathfrak{c}A}\ltimes \alpha&
    =\alpha;\label{OpActUnit1}\\
    f\ltimes\langle\textbf{id}_{\mathfrak{c}a}\rangle&
    =\textbf{id}_{\mathfrak{c}A};\label{OpActUnit2}\\
    f\cdot\sigma
    \ltimes\langle \alpha^a\rangle&
    =f\ltimes\langle \alpha^{\sigma^{-1}a'}\rangle\cdot \sigma\langle B^a\rangle;\label{OpActEquiv1}\\
    f\ltimes
    \langle \alpha^a\cdot\tau^a\rangle&
    =f\ltimes\langle \alpha^a\rangle\cdot\Pi_A\tau^a.\label{OpActEquiv2}
\end{align}
We refer to the operad $\mathcal P$ as the \textit{additive relative operad} and $\mathcal G$ as the \textit{multiplicative relative operad} of the pair.
\end{definition}

For $X=((X_d,0_d,1_d),(X_c,0_c,1_c))\in\texttt{Top}^2_{\mathds S^0}$ define
$$
    X^{\wedge-}:\mathbb S^{\text{inj}}_{\{d,c\}}\rightarrow\texttt{Top};\ X^{\wedge A}:=\wedge_A X_{\mathfrak c a},\ \sigma\cdot [x^a]:=\left[\begin{cases}
        1_{\mathfrak{c}a'},&a'\not\in\text{Im }\sigma;\\
        x^{\sigma^{-1}a'},&a'\in\text{Im }\sigma.
    \end{cases}\right]
$$
with the zeros as base points for the wedge products. We can then define the monad $(G_0;\eta,\mu)$ on $\texttt{Top}_{\mathds S^0}^2$ with
\begin{gather*}
    G_0X_\star:=\textstyle \int^{\mathbb S_\star^{\text{inj}}}\mathcal G A_+\wedge X^{\wedge A};\\
	\eta_\star x:=[\text{id}_\star,x],
	\hspace{1cm}\mu_\star[f,[[g^a,[x^{a,b}]]]]:=[f \langle g^a\rangle,[ x^{a,b}]].
\end{gather*}

\begin{definition}
	A \textit{$\mathcal G_0$-space} is a $G_0$-algebra, i.e. a pair of $\mathds S^0$-spaces $X\in\texttt{Top}^2_{\mathds S^0}$  equipped with a structural map $\chi:G_0 X\rightarrow X$ satisfying, using the notation 
	$$
	    f[ x^a]=\chi_A[f,[ x^a]]
	$$
	similar equations as in \ref{P-alg} and also that 0 is an absorbing element, ie
\begin{equation*}
     \exists a\in A: x^a=0_{\mathfrak{c} a}  \implies f[x^a]=0_{\mathfrak{c} A}.
\end{equation*}
The category of $\mathcal G_0$-spaces is denoted $\mathcal G_0[\texttt{Top}]$.
\end{definition}

If $\mathcal G$ acts on $\mathcal P$ then the functor $P$ induces a monad on $\mathcal G_0[\texttt{Top}]$.

\begin{definition}
    Let $(\mathcal G,\mathcal P)$ be a relative operad pair. A $(\mathcal G,\mathcal P)$-space is a $P$-algebra in $\mathcal G_0[\texttt{Top}]$. Equivalently a $(\mathcal G,\mathcal P)$-space is a pair of $\mathds S^0$-spaces $X\in \texttt{Top}_{\mathds S^0}^2$ equipped with a $\mathcal G_0$-space structure and a $\mathcal P$-space structure with neutral elements the zeros such that
    \begin{equation*}
        f[ \alpha^a\langle x^{a,b}\rangle]=f\ltimes\langle \alpha^a\rangle\langle f [x^{a,b^a}]\rangle.
    \end{equation*}
The category of $(\mathcal P,\mathcal G)$-spaces is denoted $(\mathcal P,\mathcal G)[\texttt{Top}]$.
\end{definition}

There is a natural operad pair structure on $(Com^\shortrightarrow,Com^\shortrightarrow)$. Set the notation $\sum_A\in Com^\shortrightarrow(A)$ for the additive copy of $Com^\shortrightarrow$ and $\prod_A\in Com^\shortrightarrow(A)$ for the multiplicative copy of  $Com^\shortrightarrow$. Then in a $(Com^\shortrightarrow,Com^\shortrightarrow)$-space the distributivity equations and the equality of the additive and multiplicative homomorphisms
\begin{align*}
    \textstyle\prod_A\textstyle\sum_{B^a}x^{a,b}&=\textstyle\sum_{\Pi_AB^a}\textstyle\prod_{A_{\langle b^a\rangle}}x^{a,b^a}\\
    \phi_+x&=\textstyle\prod_{\{ 1_c,2_c\}_c}\langle \phi_+ x , 1_c\rangle\\
    &=\textstyle\prod_{\{ 1_c,2_c\}_c}\ltimes\langle \phi_+,\textbf{id}_c\rangle\prod_{\{ 1_d,2_c\}_c}\langle x , 1_c\rangle\\
    &=\phi_\cdot x 
\end{align*}
hold. This means that $(Com^\shortrightarrow,Com^\shortrightarrow)[\texttt{Top}]$ is isomorphic to the category of topological commutative semi-algebras over commutative semi-rings\footnote{Semi-algebras and  semi-rings are like algebras and rings without the assumption that additive inverses exist, ie we have an additive commutative monoid instead of an additive abelian group}.

The main example of multiplicative relative operad we will consider is the \textit{relative linear isometries operad} $\mathscr L^{\shortrightarrow}$ with
\begin{gather*}
    \mathscr L^{\shortrightarrow}A:=\mathscr I(\oplus_A\mathds R^\infty,\mathds R^\infty);\\
    f\cdot\sigma:=\langle f_{\sigma a}\rangle,
    \hspace{1cm}\textbf{id}:=id_{\mathds R^\infty},
    \hspace{1cm}f\circ\langle  g^a\rangle:=\langle f_ag^a_b\rangle.
\end{gather*}

The classical action of the linear isometries operad on the Steiner operad induces an action on the relative versions. The extension of $\mathscr L^\shortrightarrow$ to $\mathbb S^\text{op}_{\langle d<c\rangle}$ is given by identity maps and the action maps are given by the formula
$$
    f\ltimes\langle \alpha^a \rangle:=\left[ (t,f_a\vec u^a)\mapsto f_a(\alpha^a_{b^a}(t,\vec u^a))\right].
$$

\begin{definition}
    The category of $E^\shortrightarrow_\infty$-algebras is $(\mathcal E^\shortrightarrow,\mathscr L^\shortrightarrow)[\texttt{Top}]$ for an $E_\infty^\shortrightarrow$-operad $\mathcal E^\shortrightarrow$ equipped with an action by $\mathscr L^\shortrightarrow$. .
\end{definition}

Although we give this general definition we note that there is no known non-trivial example of an $E_\infty$ operad equipped with an $\mathscr L$-action other then the Steiner operad $\mathscr H^\shortrightarrow_\infty$. Having a $q$-cofibrant, not just mixed $\Sigma$-cofibrant example would be interesting and useful, but since we can work in the mixed model structure of spectra it is not necessary.

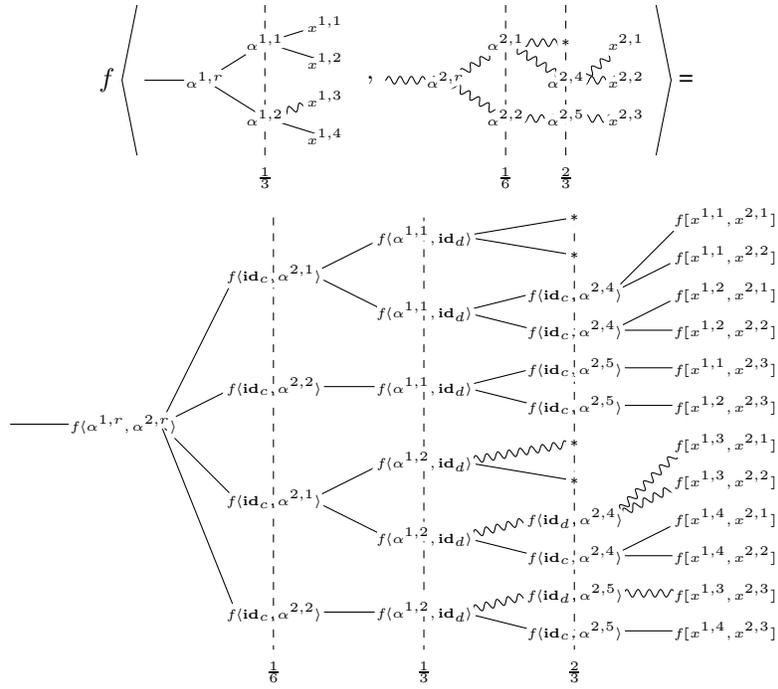
\begin{figure}
    \centering
    \begin{tikzpicture}
        
    \draw (0,0) node {$f$};
    
    \draw (0.4,-1) -- (0.2,0) --  (0.4,1);
        
    \draw (0.5,0)--(1.3,0)--(2.1,0.5)--(2.9,0.75);
    \draw (2.1,0.5)--(2.9,0.25);
    
    \draw (1.3,0)--(2.1,-0.5)--(2.9,-0.75);
    \draw[decorate,decoration={snake,segment length=1.5mm,amplitude=0.5mm}] (2.1,-0.5)--(2.9,-0.25);
        
    \draw (1.3,0) node {\tiny\contour{white}{$\alpha^{1,r}$}};
    
    \draw[dashed] (2.1,-1)--(2.1,1);
    \draw (2.1,-1.3) node {\tiny$\frac{1}{3}$};
        
    \draw (2.1,0.5) node {\tiny\contour{white}{$\alpha^{1,1}$}};
    \draw (2.1,-0.5) node {\tiny\contour{white}{$\alpha^{1,2}$}};
    
    \draw (2.9,0.75) node {\tiny\contour{white}{$x^{1,1}$}};
    \draw (2.9,0.25) node {\tiny\contour{white}{$x^{1,2}$}};
    \draw (2.9,-0.25) node {\tiny\contour{white}{$x^{1,3}$}};
    \draw (2.9,-0.75) node {\tiny\contour{white}{$x^{1,4}$}};

    \draw (3.5,0) node {$,$};
        
    \draw[decorate,decoration={snake,segment length=1.5mm,amplitude=0.5mm}] (3.7,0)--(4.5,0)--(5.3,0.5)--(6.1,0)--(6.7,0.5);
    \draw[decorate,decoration={snake,segment length=1.5mm,amplitude=0.5mm}] (5.3,0.5)--(6.1,0.5);

    \draw[decorate,decoration={snake,segment length=1.5mm,amplitude=0.5mm}] (6.1,0)--(6.9,0);
    \draw[decorate,decoration={snake,segment length=1.5mm,amplitude=0.5mm}] (4.5,0)--(5.3,-0.5)--(6.1,-0.5)--(6.9,-0.5);

    \draw (4.5,0) node {\tiny\contour{white}{$\alpha^{2,r}$}};
    
    \draw[dashed] (5.3,-1)--(5.3,1);
    \draw (5.3,-1.3) node {\tiny$\frac{1}{6}$};
    
    \draw (5.3,0.5) node {\tiny\contour{white}{$\alpha^{2,1}$}};
    \draw (5.3,-0.5) node {\tiny\contour{white}{$\alpha^{2,2}$}};
    
    \draw[dashed] (6.1,-1)--(6.1,1);
    \draw (6.1,-1.3) node {\tiny$\frac{2}{3}$};
    
    \draw (6.1,0.5) node {\tiny\contour{white}{$\ast$}};
    \draw (6.1,0) node {\tiny\contour{white}{$\alpha^{2,4}$}};
    \draw (6.1,-0.5) node {\tiny\contour{white}{$\alpha^{2,5}$}};
    
    \draw (6.9,0.5) node {\tiny\contour{white}{$x^{2,1}$}};
    \draw (6.9,0) node {\tiny\contour{white}{$x^{2,2}$}};
    \draw (6.9,-0.5) node {\tiny\contour{white}{$x^{2,3}$}};
        
    \draw (7.3,-1) -- (7.5,0) --  (7.3,1);

    \draw (7.7,0) node {$=$};
    \end{tikzpicture}
    
    \begin{tikzpicture}
    \draw (0.5,0)--(2,0);
    
    \draw (2.5,0)--(3.5,2);
    \draw (2.5,0)--(3.5,0.5);
    \draw (2.5,0)--(3.5,-1);
    \draw (2.5,0)--(3.5,-2.5);
    \draw (2,0) node {\tiny\contour{white}{$f\langle\alpha^{1,r},\alpha^{2,r}\rangle$}};

    \draw[dashed] (4,2.75)--(4,-3);
    \draw (4,-3.3) node {\tiny$\frac{1}{6}$};
    
    \draw (4.5,2)--(5.5,2.5);
    \draw (4.5,2)--(5.5,1.5);
    \draw (4,2) node {\tiny\contour{white}{$f\langle\textbf{id}_c,\alpha^{2,1}\rangle$}};
    
    \draw (4.5,0.5)--(5.5,0.5);
    \draw (4,0.5) node {\tiny\contour{white}{$f\langle\textbf{id}_c,\alpha^{2,2}\rangle$}};
    
    \draw (4.5,-1)--(5.5,-0.5);
    \draw (4.5,-1)--(5.5,-1.5);
    \draw (4,-1) node {\tiny\contour{white}{$f\langle\textbf{id}_c,\alpha^{2,1}\rangle$}};
    
    \draw (4.5,-2.5)--(5.5,-2.5);
    \draw (4,-2.5) node {\tiny\contour{white}{$f\langle\textbf{id}_c,\alpha^{2,2}\rangle$}};
    
    \draw[dashed] (6,2.75)--(6,-3);
    \draw (6,-3.3) node {\tiny$\frac{1}{3}$};
    
    \draw (6.5,2.5)--(8,2.75);
    \draw (6.5,2.5)--(8,2.25);
    \draw (6,2.5) node {\tiny\contour{white}{$f\langle\alpha^{1,1},\textbf{id}_d\rangle$}};
    
    \draw (6.5,1.5)--(7.5,1.75);
    \draw (6.5,1.5)--(7.5,1.25);
    \draw (6,1.5) node {\tiny\contour{white}{$f\langle\alpha^{1,1},\textbf{id}_d\rangle$}};
    
    \draw (6.5,0.5)--(7.5,0.75);
    \draw (6.5,0.5)--(7.5,0.25);
    \draw (6,0.5) node {\tiny\contour{white}{$f\langle\alpha^{1,1},\textbf{id}_d\rangle$}};
    
    \draw[decorate,decoration={snake,segment length=1.5mm,amplitude=0.5mm}] (6.5,-0.5)--(8,-0.25);
    \draw (6.5,-0.5)--(8,-0.75);
    \draw (6,-0.5) node {\tiny\contour{white}{$f\langle\alpha^{1,2},\textbf{id}_d\rangle$}};
    
    \draw[decorate,decoration={snake,segment length=1.5mm,amplitude=0.5mm}] (6.5,-1.5)--(7.5,-1.25);
    \draw (6.5,-1.5)--(7.5,-1.75);
    \draw (6,-1.5) node {\tiny\contour{white}{$f\langle\alpha^{1,2},\textbf{id}_d\rangle$}};
    
    \draw[decorate,decoration={snake,segment length=1.5mm,amplitude=0.5mm}] (6.5,-2.5)--(7.5,-2.25);
    \draw (6.5,-2.5)--(7.5,-2.75);
    \draw (6,-2.5) node {\tiny\contour{white}{$f\langle\alpha^{1,2},\textbf{id}_d\rangle$}};
    
    \draw[dashed] (8,2.75)--(8,-3);
    \draw (8,-3.3) node {\tiny$\frac{2}{3}$};
    
    \draw (8,2.75) node {\tiny\contour{white}{$\ast$}};
    \draw (8,2.25) node {\tiny\contour{white}{$\ast$}};
    
    \draw (8.5,1.75)--(9.5,2.75);
    \draw (8.5,1.75)--(9.5,2.25);
    \draw (8,1.75) node {\tiny\contour{white}{$f\langle\textbf{id}_c,\alpha^{2,4}\rangle$}};
    
    \draw (8.5,1.25)--(9.5,1.75);
    \draw (8.5,1.25)--(9.5,1.25);
    \draw (8,1.25) node {\tiny\contour{white}{$f\langle\textbf{id}_c,\alpha^{2,4}\rangle$}};
    
    \draw (8.5,0.75)--(9.5,0.75);
    \draw (8,0.75) node {\tiny\contour{white}{$f\langle\textbf{id}_c,\alpha^{2,5}\rangle$}};
    
    \draw (8.5,0.25)--(9.5,0.25);
    \draw (8,0.25) node {\tiny\contour{white}{$f\langle\textbf{id}_c,\alpha^{2,5}\rangle$}};
    
    \draw (8,-0.25) node {\tiny\contour{white}{$\ast$}};
    \draw (8,-0.75) node {\tiny\contour{white}{$\ast$}};
    
    \draw[decorate,decoration={snake,segment length=1.5mm,amplitude=0.5mm}] (8.5,-1.25)--(9.5,-0.25);
    \draw[decorate,decoration={snake,segment length=1.5mm,amplitude=0.5mm}] (8.5,-1.25)--(9.5,-0.75);
    \draw (8,-1.25) node {\tiny\contour{white}{$f\langle\textbf{id}_d,\alpha^{2,4}\rangle$}};
    
    \draw (8.5,-1.75)--(9.5,-1.25);
    \draw (8.5,-1.75)--(9.5,-1.75);
    \draw (8,-1.75) node {\tiny\contour{white}{$f\langle\textbf{id}_c,\alpha^{2,4}\rangle$}};
    
    \draw[decorate,decoration={snake,segment length=1.5mm,amplitude=0.5mm}] (8.5,-2.25)--(9.5,-2.25);
    \draw (8,-2.25) node {\tiny\contour{white}{$f\langle\textbf{id}_d,\alpha^{2,5}\rangle$}};
    
    \draw (8.5,-2.75)--(9.5,-2.75);
    \draw (8,-2.75) node {\tiny\contour{white}{$f\langle\textbf{id}_c,\alpha^{2,5}\rangle$}};

    \draw (10,2.75) node {\tiny\contour{white}{$f[x^{1,1},x^{2,1}]$}};
    \draw (10,2.25) node {\tiny\contour{white}{$f[x^{1,1},x^{2,2}]$}};
    \draw (10,1.75) node {\tiny\contour{white}{$f[x^{1,2},x^{2,1}]$}};
    \draw (10,1.25) node {\tiny\contour{white}{$f[x^{1,2},x^{2,2}]$}};
    \draw (10,0.75) node {\tiny\contour{white}{$f[x^{1,1},x^{2,3}]$}};
    \draw (10,0.25) node {\tiny\contour{white}{$f[x^{1,2},x^{2,3}]$}};
    \draw (10,-0.25) node {\tiny\contour{white}{$f[x^{1,3},x^{2,1}]$}};
    \draw (10,-0.75) node {\tiny\contour{white}{$f[x^{1,3},x^{2,2}]$}};
    \draw (10,-1.25) node {\tiny\contour{white}{$f[x^{1,4},x^{2,1}]$}};
    \draw (10,-1.75) node {\tiny\contour{white}{$f[x^{1,4},x^{2,2}]$}};
    \draw (10,-2.25) node {\tiny\contour{white}{$f[x^{1,3},x^{2,3}]$}};
    \draw (10,-2.75) node {\tiny\contour{white}{$f[x^{1,4},x^{2,3}]$}};
\end{tikzpicture}
\caption{$\mathscr L^\shortrightarrow$-structure of $\overline B_2X$}
\label{fig:L str of B2X}
\end{figure}

The images of $\overline B_2 X$ are also $\mathscr L^\shortrightarrow_0$-pairs with structural maps defined as
\begin{gather}\label{LActB2X}
	f[[[\alpha^{a,r},\langle \alpha^{a,v}\rangle,\langle x^{a,e}\rangle]_{T^a}, t^a]]
    :=\\
    \nonumber\left[\left[f\ltimes\langle\alpha^{a,r}\rangle,
    \left\langle f\ltimes\left\langle\begin{cases}
        \textbf{id}_{\mathfrak c s^{-1}v^a},
        &\hspace{-0.3cm}\exists i>\lvert v^a\rvert: \delta^ai= \delta^a\lvert v^a\rvert\\
        \alpha^{a,v^a},
        &\hspace{-0.3cm}\text{otherwise}
    \end{cases}\right\rangle\right\rangle
    ,\langle f[x^{a,e^a}]\rangle\right]_{\Pi_A T^a\cdot\delta^a}\hspace{-0.75cm},\lhd_A t^a\right]
\end{gather}
which is illustrated in figure \ref{fig:L str of B2X}.

\subsection{Model structure of \texorpdfstring{$E^\shortrightarrow_\infty$}{}-algebras}

The model structure of $E^\shortrightarrow_\infty$-pairs is transferred from the $q$-model structure of $\texttt{Top}_\ast^2$ by the adjunction
$$(E^\shortrightarrow L^\shortrightarrow\dashv U):\texttt{Top}_\ast^2\leftrightharpoons(\mathcal E^\shortrightarrow,\mathscr L^\shortrightarrow)[\texttt{Top}],
$$
so the weak equivalences and fibrations are respectively the maps that are $q$-equivalences and $q$-fibrations as topological space maps \cite{QuHomAlg}. All objects are fibrant and cofibrant algebras are retracts of cellular $E^\shortrightarrow_\infty$-algebras, with cells the $E^\shortrightarrow L^\shortrightarrow$-images of the cells in $\texttt{Top}_\ast^2$.

\section{Recognition of algebra spectra}

\subsection{Coordinate-free spectra}
We give a brief review of coordinate-free spectra \cite{LMS} and give some examples.

Let $\mathds U\in\mathscr I$ be countably infinite dimensional (In the context of coordinate-free spectra we refer to $\mathds U$ as a \textit{universe}). The topological category $\texttt{Sp}_{\mathds U}$ of coordinate-free $\mathds U$-spectra is composed of the class of objects
\begin{gather*}
    \{Y=(\langle Y_U\rangle,\langle\sigma^U_V\rangle)\in \Sigma_{\Pi_{\mathscr A}\texttt{Top}_\ast}\Pi_{\Sigma_{\mathscr A}\mathscr A_U}\texttt{Top}_\ast(Y_U\wedge\mathds S^{V-U},Y_{V})\mid\\
    \sigma^U_U[y,\vec 0]=y,
    \hspace{0.5cm} \sigma^V_W[\sigma^U_V[y,\vec v],\vec w]=\sigma^U_W[y,\vec v+\vec w]\}
\end{gather*}
and the morphisms spaces $\texttt{Sp}_{\mathds U}(Y,Z)$ defined as
$$
    \{\mathfrak f=\langle\mathfrak f_U\rangle\in\Pi_{\mathscr A}\texttt{Top}_\ast(Y_U,Z_U)\mid 
    \sigma^U_V[\mathfrak f_U y,\vec v]=\mathfrak f_V\sigma^U_V[y,\vec v]\}.
$$

We are particularly interested here in the case $\mathds U=\mathds R^\infty$ and in this case we use the notation $\texttt{Sp}:=\texttt{Sp}_{\mathds R^\infty}$.

\begin{example}\label{Example spectra}
    Interesting coordinate-free spectra to keep in mind are the following, with details similar to the equivalent symmetric examples in \cite[Section I.2]{Sch12}:
    \begin{itemize}
    \item For each $p\in\mathbb Z$ the $p$-sphere spectrum is defined as
    $$
    \mathds S^p:=\begin{cases}
        \langle\mathds S^{U-\mathds R^{-p}}\rangle,
        \sigma^U_V[\vec u,\vec v]:=\vec u+\vec v_{V-\mathds R^{-p}},
        & p< 0\\
        \langle\mathds S^U\rangle, \sigma^U_V[\vec u,\vec v]:=\vec u+\vec v,&p=0\\
        \langle\mathds S^{U\oplus\mathds R^p}\rangle,\ \sigma^U_V[(\vec u,\vec w),\vec v]:=(\vec u+\vec v,\vec w),& p> 0
    \end{cases}
    $$
    We use the notation $\mathds S:=\mathds S^0$.
    
    \item For each $G\in\texttt{AbGrp}$ define the \textit{Eilemberg-MacLane} spectrum
    $$
        \textstyle HG:=\langle  G\otimes F[\mathds S^U]_\ast\rangle;
        \hspace{0.5cm}\sigma^U_V[g_a\otimes \vec u^a,\vec v]:=g_a\otimes \vec u^a+\vec v.
    $$
    where $F[\mathds S^U]_\ast$ denotes the quotient of the free abelian group generated by the points of the $U$-sphere by the subgroup generated by $\infty$, and as in the Einstein convention $g_a\otimes \vec u^a$ indicates a finite sum of elements. Note that $g\otimes \infty=0$.
    
    \item For each $U\in\mathscr A$ let $O_U$ be the orthogonal group of isometric automorphisms of $U$. The total space $EO_U$ of the universal principal $O_U$-bundle is the geometric realization of the simplicial space $O_U^-\in \texttt{Top}^{\Delta^{\text{op}}}$ with
    $$
        \langle f^j\rangle\cdot\partial_i:=\left\langle\begin{cases}
            f^j,&j<i-1\\
            f^jf^{j+1},&j=i-1\\
            f^{j+1},&j> i-1
        \end{cases}\right\rangle,
        \hspace{0.5cm}
        \langle f^j\rangle\cdot\delta_i:=\left\langle\begin{cases}
            f^j,&j<i\\
            id,&j=i\\
            f^{j-1},&j>i
        \end{cases}\right\rangle.
    $$
    The $U$-spheres admit a left $O_U$-action by evaluation $f\cdot \vec u:=f\vec u$ and $EO_{U,+}$ admits the right $O_U$-action
    $$
        [\langle g^i\rangle, t]\cdot f:=\left[\left\langle\begin{cases}
            g^i,&i<m\\
            g^mf,&i=m
        \end{cases}\right\rangle, t\right].
    $$
    For $(U,V)\in\Sigma_{\mathscr A}\mathscr A_U$ we have a natural inclusion 
    $$
        \iota^U_V:O_U\rightarrow O_V,
        \hspace{1 cm} \iota^U_V f\vec v:= \vec v_{V-U}+f\vec v_U.
    $$
    We can define the \textit{Thom spectrum} as
    \begin{gather*}
        \textstyle MO:=\langle
        EO_{U,+}\wedge_{O_U}
        \mathds S^U\rangle;
        \hspace{0.25cm}\sigma^U_V[[\langle f^i\rangle, t,\vec u],\vec v]:=[\langle \iota^U_Vf^i\rangle, t,\vec u+\vec v].
    \end{gather*}
\end{itemize}
\end{example}

An $\Omega$-spectrum is a spectrum $Y\in\texttt{Sp}$ such that the adjoint structural maps $\widetilde \sigma^U_V\in\texttt{Top}_\ast(Y_U,Y_{V}^{\mathds S^{V-U}})$ are $q$-equivalences.

The stable homotopy groups of spectra are defined as $\pi_p^S Y:=\pi_0 \texttt{Sp}(\mathds S^p,Y)$. If $Y$ is an $\Omega$-spectrum then 
$$
    \pi_p^S Y\cong \begin{cases}
        \pi_0 Y_{\mathds R^{|p|}},&p<0\\
        \pi_p Y_0,& p\geq 0
    \end{cases}
$$

Spectra maps that induce isomorphisms of the stable homotopy groups are called \textit{stable weak equivalences}, and spectra $Y\in\texttt{Sp}$ with $\pi_p
^SY$ trivial for $p<0$ are called \textit{connective}.

We base space functor
\begin{equation*}
    \Lambda^\infty:\texttt{Sp}\rightarrow \texttt{Top}_\ast;\hspace{1cm}\Lambda^\infty Y:= Y_0
\end{equation*}
which is right adjoint to the suspension spectrum functor
\begin{equation*}
    \Sigma^\infty:\texttt{Top}_\ast\rightarrow \texttt{Sp};
    \hspace{0.5cm}\Sigma^\infty X:=\langle X\wedge\mathds S^{U}\rangle; 
    \hspace{0.5cm} \sigma^U_V[[x,\vec u],\vec v]:=[x,\vec u+\vec v].
\end{equation*}
with unit and counit of the adjunction
\begin{equation*}
    \eta x:=[x,\vec 0];\hspace{1cm}\epsilon_U [y,\vec u]:=\sigma^0_U[y,\vec u].
\end{equation*}

\subsection{Stable mixed model structure of spectra} For any spectrum $Y\in\texttt{Sp}$ the cylinder spectrum is defined as $Y\wedge I_+:=\langle Y_U\wedge I_+\rangle$ and the cone spectrum is $CY:=Y\wedge I_+/_{[y,1]\sim [y',1]}$.

In the strict Quillen model structure on $\texttt{Sp}$ a morphism $\mathfrak f\in\texttt{Sp}(X,Y)$ is a weak equivalence if each $\mathfrak f_U$ is a $q$-equivalence, a fibration if it is a Serre fibration, ie if it has the homotopy lifting property with respect to the cylinder inclusions of cones of sphere
spectra $\text{in}_0\in\texttt{Sp}(C\mathds S^q,C\mathds S^q\wedge I_+)$ for all $q\in\mathds Z$, and a cofibration if it is a retract of a relative cell-spectrum, with cells given by cones of sphere spectra and domain of the attaching maps the boundary sphere spectra \cite[Section VII.4]{EKMM}. This is a cofibrantly generated model structure with factorization systems induced by the small object argument. The weak equivalences, fibrations and cofibrations of this model structure are referred to as $q$-equivalences, $q$-fibrations and $q$-cofibrations respectively.

Homotopy equivalences in $\texttt{Sp}$ are spectra maps that admit an inverse up to homotopy, with homotopies defined via the cylinder spectra in the usual way. In the strict Hurewicz/Str\o m model structure $\mathfrak f$ is a weak equivalence if it is a homotopy equivalence, a fibration if it is a Hurewicz fibration, ie if it has the homotopy lifting property with respect to all cylinder inclusions $\text{in}_0\in\texttt{Sp}(X,X\wedge I_+)$, and a cofibration if it has the left lifting property against trivial Hurewicz fibrations.

The weak factorization system can be constructed through (co)monads as described in \cite{BR13}. For any $Y\in\texttt{Sp}$ let the spectrum of Moore paths in $Y$ be
\begin{gather*}
    MY:=\langle \Sigma_{[0,\infty)}\{Y_U^{[0,\infty]_+}\mid s\geq t\implies \gamma s=\gamma t\}\rangle,\\
    \sigma^U_V[(t,\gamma),\vec v]:= (t,r\mapsto \sigma^U_V [\gamma r,\vec v]).
\end{gather*}
The factorization systems are then defined as
\begin{gather*}
    \Gamma\mathfrak f:=X\times_Y MY; 
    \hspace{0.5cm}C_t\mathfrak f x:=(x,0,r\mapsto \mathfrak f x),
    \hspace{0.5cm}F\mathfrak f(x,t,\gamma):=\gamma t.\\
    E\mathfrak f:=\Gamma\mathfrak f\wedge[0,\infty]_+\sqcup_{\Gamma \mathfrak f}Y;
    \hspace{0.5cm}C\mathfrak f x:=(x,0,r\mapsto \mathfrak f x,0),
    \hspace{0.5cm}F\mathfrak f(x,t,\gamma,s):=\gamma s.
\end{gather*}

The weak equivalences, fibrations and cofibrations of this model structure are referred to as $h$-equivalences, $h$-fibrations and $h$-cofibrations respectively. We then equip $\texttt{Sp}$ with the mixed model structure as described in \cite[Prop 3.6]{Cole06}.

Since the point of spectra is to study stabilization phenomena we are actually interested in inverting the stable weak homotopy equivalences. The stable model structure with stable weak homotopy equivalences as weak equivalences is obtained from the strict model structure by the process of Bousfield localization through the following idempotent monad \cite{BF78,Sc97}. For every spectrum $Y\in\texttt{Sp}$ we can functorialy define an inclusion spectrum\footnote{Inclusion spectra are those with adjoint structural maps $\tilde \sigma$ all inclusions.} $\widetilde Y$  equipped with a quotient map $Y\rightarrow \widetilde Y$, so we may think of points in $\widetilde Y$ as equivalence classes of points in $Y$ (see \cite[Ap1]{LMS} for a detailed construction). If $Y$ is already an inclusion spectrum then $\widetilde Y=Y$. We may then define the \textit{spectrification functor}
\begin{equation*}
    \widetilde \Omega:\texttt{Sp}\rightarrow\texttt{Sp}; \hspace{0.3cm}
    \widetilde\Omega Y:=\langle\text{colim}_{\mathscr A_U}\widetilde Y_{V}^{\mathds S^{V-U}}\rangle; \hspace{0.3cm}
    \sigma^U_W[\gamma,\vec w]:=[\vec v\mapsto \gamma(\vec v+\vec w)]
\end{equation*}
induced by the adjoint structural maps $\tilde \sigma$ and with the formula for the structural maps determined for a choice of representative $\gamma$ with domain $V\in\mathscr A_W$. This is a Quillen idempotent monad with structural natural map 
\begin{equation}\label{Stabilization monad map}
    \epsilon':Id\Rightarrow \widetilde \Omega;
    \hspace{1cm}\epsilon_U' y:=\left[\vec v\mapsto \sigma^U_V[y,\vec v]\right] 
\end{equation}

The stable model structure on spectra $\texttt{Sp}_{\widetilde \Omega}$ has as weak equivalences the stable weak equivalences and stable fibrations are $\mathfrak p\in\texttt{Sp}(E,B)$ composed of indexwise Hurewicz fibrations such that the maps
$$
    (\tilde\sigma^U_V,\mathfrak p_U):E_U\rightarrow E_{V}^{\mathds S^{V-U}}\times_{B_{V}^{\mathds S^{V-U}}}B_{U}
$$
are $q$-equivalences. The fibrant spectra are the $\Omega$-spectra, and the cofibrant spectra are those homotopy equivalent to retracts of $q$-cofibrant spectra.  With the induced stable model structure the adjunction $(\Sigma^\infty\dashv \Lambda^\infty)$ is a Quillen adjunction.

The morphisms category $\texttt{Sp}^\shortrightarrow$ admits a projective stable model structure with $(\mathfrak f_d,\mathfrak f_c)\in \texttt{Sp}^\shortrightarrow(\mathfrak i:Y_d\rightarrow Y_c,\mathfrak j:Z_d\rightarrow Z_c)$ a weak equivalence or fibration if $\mathfrak f_c$ and $\mathfrak f_o$ are both stable weak equivalences or stable fibrations respectively, and it is a cofibration if both $\mathfrak f_d$ and $(\mathfrak f_c,\mathfrak j):Y_c\vee_{Y_d}Z_d\rightarrow Z_c$ are stable cofibrations.

\subsection{Recognition of \texorpdfstring{$\infty$}{}-loop maps} We can now prove the recognition principle for $\infty$-loop pairs of spaces of spectra maps. The base pair of spaces functor is
\begin{equation*}
    \Lambda_2^\infty:\texttt{Sp}^\shortrightarrow\rightarrow\texttt{Top}_\ast^2, 
    \hspace{0.5cm}\Lambda_2^\infty\mathfrak i:=(Y_{d,0},Y_{c,0})
\end{equation*} 
and the relative suspension functor is
\begin{equation*}
    \Sigma^\infty_\shortrightarrow:\texttt{Top}_\ast^2\rightarrow\texttt{Sp}^\shortrightarrow,
    \hspace{0.5cm}\Sigma^\infty_\shortrightarrow X:= \Sigma^\infty (\text{in}_d: X_d\rightarrow X_d\vee X_c).
\end{equation*}
 
We have a Quillen adjunction 
\begin{equation*}
    (\Sigma^\infty_\shortrightarrow\dashv\Lambda^\infty_2):\texttt{Top}_\ast^2\leftrightharpoons \texttt{Sp}^\shortrightarrow;\ \eta_\star x:=[x,\vec 0],\ \epsilon_{\star,U} [y,\vec u]:=\left[\begin{cases}
        \sigma^0_U[y,\vec u],&\mathfrak{c}y=\star\\
        \sigma^0_U[\mathfrak i y,\vec u],&\mathfrak{c}y\neq \star
    \end{cases}\right].
\end{equation*}

The spectrification functor $\widetilde \Omega$ induces
\begin{equation*}
    \widetilde \Omega_\shortrightarrow:\texttt{Sp}^\shortrightarrow\rightarrow\texttt{Sp}^\shortrightarrow;
    \hspace{0.5cm}\widetilde \Omega_\shortrightarrow\mathfrak i:=(\widetilde \Omega\mathfrak i:\widetilde \Omega Y_d\rightarrow \widetilde \Omega Y_c).
\end{equation*}
The $\infty$-loop pair of spaces functor is defined as
\begin{equation*}
    \Omega_2^\infty:\texttt{Sp}^\shortrightarrow\rightarrow \mathcal E^\shortrightarrow[\texttt{Top}];
    \hspace{0.5cm}\Omega_2^\infty\mathfrak i:=\Lambda_2^\infty\widetilde\Omega_\shortrightarrow\mathfrak i
\end{equation*}
with structural maps induced by the formula \ref{SumRelLopSpc} by taking representatives of the $\gamma^a$ with a common domain.

This functor is not a right adjoint, but it is a weak Quillen right quasiadjoint. The left quasiadjoint functor is defined as follows: We have simplicial pointed maps $B_-(\Sigma^U_\shortrightarrow,E_U^\shortrightarrow,X)\in(\texttt{Top}_\ast^\shortrightarrow)^\Delta$ with
\begin{gather*}
    \textstyle B_{\langle m\rangle}(\Sigma^U_\shortrightarrow,E_U^\shortrightarrow,X)_\star\cong
    (\int^{\mathbb T^0_\star \langle m\rangle}\Pi_V\mathcal E^\shortrightarrow\text{in }v\times \Pi_{E^m}X_{\mathfrak ce})\wedge \mathds S^U,\\
    [\langle \alpha^v\rangle,\langle x^e\rangle,\vec u]_T\cdot \partial_i:=
    \begin{cases}
        \begin{cases}
            [\langle \alpha^v\rangle,\langle x^e\rangle,\alpha^{1,-1}_{e'}\vec u]_{T_{\geq e'}},&\vec u\in\alpha^1_{e'} U\\
            \infty,&\vec u\not\in\alpha^1\sqcup_{E^0} U
        \end{cases},&i=0\\
        \left[\left\langle \begin{cases}
            \alpha^v\langle \alpha^{se'}\rangle,&\lvert v\rvert=i-1\\
            \alpha^v,&\lvert v\rvert\neq i-1
        \end{cases} \right\rangle,\langle x^e\rangle,\vec u\right]_{T\cdot\partial_i},&0<i<m\\
        [\langle \alpha^v\rangle,\langle\alpha^{s^{-1}e}\langle x^{e'}\rangle\rangle,\vec u]_{T\cdot\partial_m}&,i=m
    \end{cases}\\
    [\langle \alpha^v\rangle,\langle x^e\rangle,\vec u]_T\cdot \delta_i:=\left[\left\langle \begin{cases}
        \textbf{id}_{\mathfrak c s^{-1}v},&\lvert v\rvert= i\\
        \alpha^v,&\lvert v\rvert\neq i
    \end{cases}\right\rangle,\langle x^e\rangle,\vec u\right]_{T\cdot \delta_i}.
\end{gather*}

Define the relative $\infty$-delooping functor as
\begin{gather*}
    \textstyle B^\infty_\shortrightarrow:\mathcal E^\shortrightarrow[\texttt{Top}]\rightarrow \texttt{Sp}^\shortrightarrow;
    \hspace{1cm} B^\infty_\shortrightarrow X_\star
    :=\langle \lvert B_-(\Sigma^U_\shortrightarrow,E_U^\shortrightarrow, X)_\star\rvert\rangle,\\
    \sigma^U_V[[[\langle\alpha^v\rangle,\langle x^e\rangle,\vec u]_T, t],\vec v]:=[[\langle\alpha^v\rangle,\langle x^e\rangle,\vec u+\vec v]_T, t].
\end{gather*}

Points in $B^\infty_\shortrightarrow X_{\star,U}$ are equivalence classes of decorated filtered rooted relative trees as in the description of the bar resolution $\overline B_2 X$, except the root vertex is decorated with a vector in $U$ and the relative operad points decorating the inner vertices must be contained in the suboperad $\mathcal E^\shortrightarrow_U$ of the $\mathscr A$-filtration of $\mathcal E^\shortrightarrow$.

\begin{theorem}\label{MorphRecogPrin}
    For $\mathcal E^\shortrightarrow$ an $E^\shortrightarrow_\infty$-operad there is an idempotent quasiadjunction
    $$
    (B_\shortrightarrow^\infty\dashv_{\ \overline B_2,\widetilde\Omega_\shortrightarrow} \Omega^\infty_2):\mathcal E^\shortrightarrow[\texttt{Top}]\leftrightharpoons\texttt{Sp}^\shortrightarrow
    $$
\end{theorem}

\textbf{Proof:} The unit span and cospan has $\eta'$ the natural weak equivalence \ref{eta'}, $\epsilon'$ induced by the idempotent monad transformation \ref{Stabilization monad map} and $\eta$ and $\epsilon$ are defined by the following formulas:
\begin{gather*}
     \eta:\overline B \Rightarrow\Omega^\infty_2B^\infty_\shortrightarrow,
     \hspace{1cm} \epsilon:B^\infty_\shortrightarrow\Omega^\infty_2\Rightarrow
    \widetilde \Omega_\shortrightarrow,\\
    \eta_{\star} [[\alpha^r,\langle \alpha^v\rangle,\langle x^e\rangle]_T, t]:=
    \left[\vec u\mapsto \begin{cases}
    	[[\langle \alpha^v\rangle,\langle x^e\rangle,\alpha_{e'}^{r,-1}\vec u]_{T_{\geq e'}}, t],
        &\vec u\in \alpha^r_{e'}U\\
        \infty,
    	&\vec u\not\in \alpha^r\sqcup_{E^0}U
        \end{cases}\right],\\
    \epsilon_{\star,U}[[\langle \alpha^v\rangle,\langle \gamma^e\rangle,\vec u]_T, t]
    :=\left[\vec v\mapsto \circ_T\alpha^v\langle\gamma^e\rangle(\vec u+\vec v)\right].
\end{gather*}

\begin{figure}
    $$\begin{tikzpicture}
        \draw (-1.4,0.05) node {$\vec u\mapsto$};
        
        \draw[decoration={brace},decorate] (-1,-1.9) --  (-1,1.9);
    
        \draw (-0.75,1.5)--(0,1.5);
        \draw [decorate,decoration={snake,segment length=1.5mm,amplitude=0.5mm}] (0,1.5) -- (1,1.5) -- (2,1.5) -- (3,1.5)--(4,1.5);
        \draw[decorate,decoration={snake,segment length=1.5mm,amplitude=0.5mm}] (3,1.5)--(4,2);
        \draw[decorate,decoration={snake,segment length=1.5mm,amplitude=0.5mm}] (2,1.5) -- (3,1)--(4,1);
        
        \draw [decorate,decoration={snake,segment length=1.5mm,amplitude=0.5mm}]
        (1,0)--(2,0.5)--(3,0.5)--(4,0.5);
        \draw (-0.75,0)--(1,0)--(2,0)--(3,0)--(4,0);
        \draw[decorate,decoration={snake,segment length=1.5mm,amplitude=0.5mm}] (3,0)--(4,-0.5);
        
        \draw (-0.75,-1.5) -- (1,-1.5) -- (2,-1.5) -- (3,-1.5)--(4,-2);
        \draw[decorate,decoration={snake,segment length=1.5mm,amplitude=0.5mm}] (3,-1.5)--(4,-1.5);
        \draw (2,-1.5) -- (3,-1)--(4,-1);
        
        \draw[dashed] (1,1.7)--(1,-2.4);
        \draw (1,-2.5) node {$t^0$};
        
        \draw (1.02,1.52) node {\contour{white}{$\alpha^1$}};
        \draw (1.02,0.02) node {\contour{white}{$\alpha^2$}};
        \draw (1.02,-1.43) node {\contour{white}{$\alpha^3$}};
        
        \draw[dashed] (2,1.7)--(2,-2.4);
        \draw (2,-2.5) node {$t^1$};
        
        \draw (2.02,1.52) node {\contour{white}{$\alpha^4$}};
        \draw (2.02,0.52) node {\contour{white}{$\alpha^5$}};
        \draw (2.02,0.02) node {\contour{white}{$\alpha^6$}};
        \draw (2.02,-1.43) node {\contour{white}{$\alpha^7$}};
        
        \draw[dashed] (3,1.7)--(3,-2.4);
        \draw (3,-2.5) node {$t^2$};
        
        \draw (3.02,1.52) node {\contour{white}{$\alpha^8$}};
        \draw (3.02,1.02) node {\contour{white}{$\alpha^9$}};
        \draw (3.02,0.52) node {\contour{white}{$\alpha^{10}$}};
        \draw (3.02,0.02) node {\contour{white}{$\alpha^{11}$}};
        \draw (3.02,-0.98) node {\contour{white}{$\alpha^{12}$}};
        \draw (3.02,-1.48) node {\contour{white}{$\alpha^{13}$}};
        
        \draw (4.45,2.15) node {$x^1$};
        \draw (4.45,1.65) node {$x^2$};
        \draw (4.45,1.15) node {$x^3$};
        \draw (4.45,0.65) node {$x^4$};
        \draw (4.45,0.15) node {$x^5$};
        \draw (4.45,-0.35) node {$x^6$};
        \draw (4.45,-0.85) node {$x^7$};
        \draw (4.45,-1.35) node {$x^8$};
        \draw (4.45,-1.85) node {$x^9$};
        
        \draw (0,1.5) node {\contour{white}{$\alpha_1^{r,-1}\vec u$}};
        
        \draw (0,0) node {\contour{white}{$\alpha_2^{r,-1}\vec u$}};
        
        \draw (0,-1.5) node {\contour{white}{$\alpha_3^{r,-1}\vec u$}};
        
        \draw[densely dotted] (-0.85,0.85)--(6.3,0.85);
        \draw[densely dotted] (-0.85,-0.65)--(6.3,-0.65);
    \end{tikzpicture}$$
        \caption{Representative $U$-loop of $\eta_c [[\alpha^r,\langle \alpha^v\rangle,\langle x^e\rangle]_T, t]$}
\end{figure}
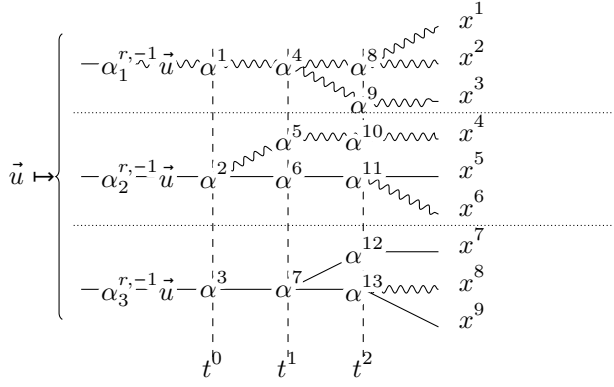

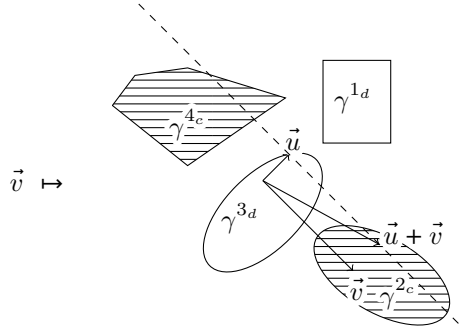
\begin{figure}
    $$\begin{tikzpicture}
        \draw (-3,0) node {\contour{white}{$\vec v\ \mapsto $}};
        
        \draw (0.8,0.5) rectangle (1.7,1.6);
        \draw[pattern=horizontal lines] (-2,1)--(-1,0.2)--(0.3,1.1)--(-1,1.5)--(-1.7,1.4)--(-2,1);
        \draw[pattern=horizontal lines,rotate=-30] (2,-0.3) ellipse (1cm and 0.5cm);
        \draw[rotate=45] (-0.3,-0.3) ellipse (1cm and 0.5cm);
        
        \draw[->] (0,0)--(0.35,0.35);
        \draw (0.4,0.55) node {\contour{white}{$\vec u$}};
        
        \draw[->] (0,0)--(1.2,-1.2);
        \draw (1.25,-1.45) node {\contour{white}{$\vec v$}};
        
        \draw[->] (0,0)--(1.55,-0.85);
        \draw (2,-0.7) node {\contour{white}{$\vec u+\vec v$}};
        
        \draw[dashed] (-1.65,2.35)--(2.6,-1.9);
        
        \draw (1.2,1.1) node {\contour{white}{$\gamma^{1_d}$}};
        
        \draw (1.8,-1.5) node {\contour{white}{$\gamma^{2_c}$}};
        
        \draw (-0.3,-0.5) node {\contour{white}{$\gamma^{3_d}$}};
        
        \draw (-1,0.75) node {\contour{white}{$\gamma^{4_c}$}};
    \end{tikzpicture}
    $$
    \caption{Representative $V$-loop of $\epsilon_{c,U}[[\langle\alpha^v\rangle,\langle \gamma^e\rangle,\vec u]_T, t]$}
    \label{fig:Represent Eta}
\end{figure}

We verify that the conditions for definition \ref{QsiAdjQuilFr} are satisfied.

(i): By the assumptions on $\mathcal E^\shortrightarrow$ and \cite[Prop. 3.2.3]{Vi20} the functor $B^\infty_\shortrightarrow$ is left derivable.

(ii): Trivially $\Omega^\infty_2$ preserves fibrant objects. Since $\Omega^\infty=\Lambda^\infty\widetilde \Omega$ and stable weak equivalences are by definition maps whose images under $\widetilde\Omega$ are strict weak equivalences we have that $\Omega^\infty_2$ preserves weak equivalences.

(iii): The functor $\overline B_2$ preserves cofibrant objects by \cite[Prop. 3.2.3]{Vi20} and trivially preserves fibrant objects. The functor $\widetilde\Omega_\shortrightarrow$ preserves cofibrant objects by the results in \cite[Sec. 5.3]{FP90} and the fact that we are using the mixed stable model structure on spectra, and it trivially preserves fibrant objects since it is the fibrant replacement functor of the stable model structure.

(iv): As a map of topological spaces $\eta'$ is a realization of a simplicial strong deformation retract, so it is itself a strong deformation retract of topological spaces and therefore in particular a $q$-equivalence \cite[theorems 9.10, 9.11 and 11.10]{MayGeomItLoopSp}. The map $\epsilon'$ is a weak equivalence by the definition of the stable model structure.

(v): The natural homotopy $H$ which gives the homotopy commutativity in $\texttt{Sp}^\shortrightarrow$
\begin{gather*}
    \epsilon_{B^\infty_\shortrightarrow X}B^\infty_\shortrightarrow\eta_X[[\langle\alpha^v\rangle,\langle[\beta^{e,r},\langle\beta^{e,w}\rangle,\langle x^{e,f}\rangle]_{S^e}, s^e]\rangle,\vec u]_T, r]\\
    =\left[\vec v\mapsto
    \begin{cases}
    	[[
        \langle \beta^{e,w}\rangle,\langle x^{e,f}\rangle,(\circ_{<e}\alpha^{v}_{e'}\beta^{e,r}_{f'})^{-1} \vec u+\vec v ]_{ S^e_{\geq f'}}, s^e],&\vec u+\vec v\in \circ_{<e}\alpha^{v}_{e'}\beta^{e,r}_{f'}V\\
        \infty,&\hspace{-1.65cm}\vec u+\vec v\not\in \circ_T\alpha^v\langle\beta^{e,r}\rangle\sqcup_{\Sigma_{E^m}F^{e,0}}V
    \end{cases}\right]
    \\
    \simeq_{H_X}
    [ \vec v\mapsto
    [[\langle\alpha^v\rangle,\langle\circ_{S^e}\beta^{e,w}\langle x^{e,f}\rangle\rangle,\vec u+\vec v
    ]_{S^e}, r]
    ]\\
    =\epsilon'_{B^\infty_\shortrightarrow X}B^\infty_\shortrightarrow\eta'_X[[\langle\alpha^v\rangle,\langle[\beta^{e,r},\langle\beta^{e,w}\rangle,\langle x^{e,f}\rangle]_{S^e},  s^e]\rangle,\vec u]_T, r]
\end{gather*}
is
\begin{gather*}
    H:B^\infty_\shortrightarrow\overline B_2 \wedge I_+\Rightarrow \widetilde \Omega B^\infty_\shortrightarrow ,\\
    H_{X,U}([[\langle\alpha^v\rangle,\langle[\beta^{e,r},\langle\beta^{e,w}\rangle,\langle x^{e,f}\rangle]_{S^e}, s^e]\rangle,\vec u]_T, r],t):=\\
    \left[\vec v\mapsto \left[\left[\left\langle\begin{cases}
        \alpha^v\\
        \textbf{id}\\
        \beta^{e,w}
    \end{cases}
    \right\rangle,\langle x^{e,f}\rangle,\vec u+\vec v\right]_{T\circ\langle S^e\cdot\delta^e\rangle} ,
    \Phi t \right]\right]
\end{gather*}
where
$$
    \Phi t:=(1-t)(\circ_{j=m+1}^{m+n+1}\partial_j\cdot r)^i+t(\circ_{j=0}^{m}\partial_0\cdot\lhd_{E^m} s^e)^i,
$$
with the conditions in the formula similar to the ones in \ref{EstructBarB}. 

(vi): In $\mathcal E^\shortrightarrow[\texttt{Top}]$ we have strict commutativity
\begin{align*}
    \Omega^\infty_2\epsilon \eta_{\Omega^\infty_2}[[\alpha^r,\langle\alpha^v\rangle,\langle \gamma^e\rangle]_T, t]&=
    [\vec v\mapsto \circ_T\alpha^v\langle\gamma^e\rangle \vec v]\\
    &=\Omega^\infty_2\epsilon' \eta'_{\Omega^\infty_2}[[\alpha^r,\langle\alpha^v\rangle,\langle \gamma^e\rangle]_T, t].\blacksquare
\end{align*}

\begin{theorem}\label{Idempotent}
    The quasiadjunction in theorem \ref{MorphRecogPrin} is idempotent and induces an equivalence
    $$
        (\mathbb LB^\infty_\shortrightarrow\dashv \mathbb R\Omega^\infty_2):\mathcal Ho \mathcal E^\shortrightarrow[\texttt{Top}]_{\text{grp}}\leftrightharpoons\mathcal Ho \texttt{Sp}^\shortrightarrow_{\text{con}}.
    $$
    between the homotopy categories of grouplike $\mathcal E^\shortrightarrow$-pairs and maps between connective spectra.
\end{theorem}

\textbf{Proof:} In $\mathcal E^\shortrightarrow[\texttt{Top}]$ the conditions for definition \ref{QsiMonIdempQuil} are satisfied and the resulting reflective homotopy subcategory is composed of the grouplike $\mathcal E^\shortrightarrow$-pairs:

(i) As we have seen $\eta'$ is a natural weak equivalence and by definition $\text{cof}$ is a natural trivial fibration, so $\text{cof}\eta'_{\mathfrak C}$ is a weak equivalence.
    
(ii) Since $\Omega^\infty_2$ preserves weak equivalence between fibrant objects and $B^\infty_\shortrightarrow$ preserves weak equivalences between cofibrant objects we have that $\Omega^\infty_2\mathfrak FB^\infty_\shortrightarrow\mathfrak C$ preserves weak equivalences.
    
(iii) The natural transformation $\eta$ is a natural group completion, since it is a realization of a simplicial group completion map (see \cite[Theorems 2.7, 9.10 and 9.11]{MayGeomItLoopSp}, and \cite[Theorem 2.2]{MayEspcGrpComplPermCat}), and the images of $\Omega^\infty_2\mathfrak FB^\infty_\shortrightarrow\mathfrak C$ are grouplike, therefore $\eta_{\Omega^\infty_2\mathfrak FB^\infty_\shortrightarrow\mathfrak C}$ is a natural weak equivalence. By naturality $\Omega^\infty_2\mathfrak FB^\infty_\shortrightarrow\mathfrak C\eta$ is also a group completion, and since the domain and codomain are grouplike this is a natural weak equivalence.

(iv) This condition holds since fibrations are preserved by pullbacks, fibrations induce long exact sequences of homotopy groups and for a fibration $p:E\twoheadrightarrow B$ and a map $f:X\rightarrow B$ the fibers of the pullback $f^\ast p:X\times_BE\rightarrow X$ are homeomorphic to the fibers of $p$.
    
(v) Pushouts in $\mathcal E^\shortrightarrow[\texttt{Top}]$ by a cofibration whose domain is $m$-cofibrant in $\texttt{Top}_\ast$ is a retract of a transfinite composition of pushouts by $m$-cofibrations in $\texttt{Top}_\ast$ (see \cite[I.4]{Sp04}), hence this condition holds since $\texttt{Top}_\ast$ with the mixed model structure is left proper and the underlying functor of $\mathcal E^\shortrightarrow$ is an $m$-cofibrant $\mathbb S_{\{d,c\}}$-space.

By the characterization of fibrations in the resulting Bousfield localization in \cite[Prop. 2.3.6]{Vi20} the fibrations are the group completions and fibrant objects are the grouplike $\mathcal E^\shortrightarrow$-pairs.

The dual conditions for definition \ref{QsiMonIdempQuil} are also satisfied in $\texttt{Sp}^\shortrightarrow$ and the resulting coreflective homotopy subcategory is composed of the maps between connective spectra. Note that conditions (i), (ii) and (iii) are self dual.

(i) By definition of the stable model structure $\epsilon'$ is a natural stable weak equivalence and by definition $\text{fib}$ is a natural trivial cofibration, so $\eta'_{\mathfrak F}\text{fib}$ is a weak equivalence.
    
(ii) That $B^\infty_\shortrightarrow\mathfrak C\Omega^\infty_2\mathfrak F$ preserves weak equivalences follows by the same argument for $\Omega^\infty_2\mathfrak FB^\infty_\shortrightarrow\mathfrak C$.
    
(iii) We have that $\eta_{\Omega^\infty_2}$ is a natural weak equivalence, and since $\Omega^\infty_2\epsilon\eta_{\Omega^\infty_2}=\Omega^\infty_2\epsilon'\eta'_{\Omega^\infty_2}$ and $\Omega^\infty_2\epsilon'\eta'_{\Omega^\infty_2}$ is a natural weak equivalence by the 2-out-of-3 property $\Omega^\infty_2\epsilon$ is a natural weak equivalence. Since the images of $\widetilde \Omega$ are $\Omega$-spectra by the formula for stable homotopy groups of $\Omega$-spectra we have that $\epsilon$ induces isomorphisms on the non-negative stable homotopy groups, and is therefore a stable weak equivalence on the maps between connective spectra. The images of $B^\infty_\shortrightarrow$ are connective by \cite[11.12]{MayGeomItLoopSp} and \cite[A5]{MayEspcGrpComplPermCat}. Therefore $\epsilon_{B^\infty_\shortrightarrow\mathfrak C\Omega^\infty_2\mathfrak F}$ is a natural weak equivalence. By naturality $B^\infty_\shortrightarrow\mathfrak C\Omega^\infty_2\mathfrak F\epsilon$ also induces isomorphisms on the non-negative stable homotopy groups and so is also a natural weak equivalence.

(iv) This condition holds since cofibrations are preserved by pullbacks, spectra cofibrations induce long exact sequences of stable homotopy groups and for any cofibration $\mathfrak i:A\hookrightarrow X$ and map $\mathfrak f:A\rightarrow Y$ the cofiber of the pushout $\mathfrak f_\ast\mathfrak i:Y\rightarrow X\sqcup_AY$ are homeomorphic to the cofibers of $\mathfrak i$.
    
(v) The stable model structure of spectra is right proper so the dual of (v) holds.

By the dual of the characterization in \cite[Prop. 2.3.6]{Vi20} the cofibrant objects are the spectra maps such that
\begin{equation*}
    \Gamma((\epsilon B^\infty_\shortrightarrow\text{cof}_{\Omega^\infty_2})_{\mathfrak F}\mathfrak i)\times_{\widetilde\Omega\mathfrak F\mathfrak i}\mathfrak i\rightarrow \mathfrak i
\end{equation*}
are weak equivalences, which is equivalent to $\iota$ being a map of connective spectra.$\blacksquare$

\subsection{\texorpdfstring{$\mathds S$}{}-modules and commutative algebra spectra}

In order to define a monoidal category of spectra, so that we get natural definitions of spectral algebraic structures, we need to work on the more structured category of sphere modules $\texttt{Mod}_{\mathds S}$ \cite{EKMM}. As a first step consider for $A\in\mathbb S$ the \textit{external smash product} functor
\begin{gather*}
    \overline \wedge_A:\Pi_A\texttt{Sp}\rightarrow \texttt{Sp}_{\oplus_A\mathds R^\infty},
    \hspace{1cm}
    \overline \wedge_A\langle Y^a\rangle:=\left\langle\wedge_AY^a_{U^a}\right\rangle,\\
    \sigma^{\langle U^a\rangle}_{\langle V^a\rangle}[[y^a],\langle \vec v^a\rangle]:=[\sigma^{U^a}_{V^a}[y^a,\vec v^a]].
\end{gather*}

The change of universe in this product is formally problematic, and the following construction is used to internalize the smash product in $\texttt{Sp}$. For $K\subset_{\text{cpct}}\mathscr L A$ define the monotone functions
\begin{align*}
    \mu\in\texttt{POSet}(\mathscr A_{\oplus_A\mathds R^\infty},\mathscr A),&\ \  
    \mu\langle U^a\rangle:=\textstyle\sum_K f\langle U^a\rangle\\
    \nu\in\texttt{POSet}(\mathscr A,\mathscr A_{\oplus_A\mathds R^\infty}),&\ \ 
    \nu  U:= \cap_K f^{-1}U
\end{align*}
which satisfy
\begin{equation*}
    \mu \nu  U\subset U,\ \ 
    \nu\mu\nu U=\nu U,\ \ 
    \langle U^a\rangle \subset \nu \mu  \langle U^a\rangle,\ \ 
    \mu\langle U^a\rangle =\mu \nu \mu  \langle U^a\rangle.
\end{equation*}

For all $(\langle U^a\rangle,V)\in \Sigma_{\mathscr A_{\oplus_A\mathds R^\infty}}\mathscr A_{\mu \langle U^a\rangle}$ we have the associated Thom complex
$$
    TK_V^{\langle U^a\rangle}:=\Sigma_K \mathds S^{V- f\langle U^a\rangle}/_{(f,\infty)\sim (g,\infty)}\in\texttt{Top}_\ast,
$$
where $\Sigma_K \mathds S^{V- f\langle U^a\rangle}$ is topologized as a subspace of $K\times\mathds S^V$, with the equivalence class $[f,\infty]$ as base point. We will use the notation ${}_{f}^{\vec v}:=[f,\vec v]\in TK^{\langle U^a\rangle}_{V}$.

The \textit{twisted half-smash product} is defined as
\begin{gather*}
    \mathscr L A \ltimes-:\texttt{Sp}_{\oplus_A\mathds R^\infty}\rightarrow \texttt{Sp};
    \hspace{0.4cm}\mathscr L A \ltimes Z:=
    \left\langle\underset{K\subset_{cpct}\mathscr L A }{\text{Colim }} TK^{\nu U}_U\wedge Z_{\nu  U}\right\rangle,\\
    \sigma^U_V [ [\prescript{\vec u}{f}{},z ],\vec v ]:= [\prescript{(\vec u+\vec v)_{V-f\nu  V}}{f}{},\sigma^{\nu U}_{\nu V} [z,f^{-1}(\vec u+\vec v)_{f\nu  V}]].
\end{gather*}

We define the monad $(\mathds L;\eta,\mu)$ on $\texttt{Sp}$ with
\begin{equation*} 
    \mathds L Y:=\mathscr{L}\underline 1 \ltimes Y;
    \hspace{0.5cm}\eta y:=[\prescript{\vec 0}{id}{},y],
    \hspace{0.5cm}\mu[\prescript{\vec u}{f}{},[\prescript{\vec v}{g}{},y]]:=[\prescript{\vec u+f\vec v}{fg}{},y].
\end{equation*}

We refer to the $\mathds L$-algebras as $\mathds L$-spectra and for $(Y,\mathfrak y)\in\mathds L[\texttt{Sp}]$ we use the notation $\prescript{\vec u}{f}{y}:=\mathfrak y[{}_f^{\vec u},y]$. 

The sphere spectrum $\mathds S$, the Eilenberg-Maclane spectra $HG$ and the Thom spectrum $MO$ in example \ref{Example spectra}, as well as the suspensions $\Sigma^\infty X$ of $\mathscr L$-spaces, deloopings $B^\infty X$ of $E_\infty$-rings and the spectrifications $\widetilde\Omega Y$ of $\mathds L$-spectra, are all $\mathds L$-spectra with structural morphisms given respectively by:
$$
    \begin{array}{|c|c|}\hline
         & \prescript{\vec u}{f}{y} \\\hline
         \mathds S&\vec u+f\vec v\\
         HG&g_a\otimes \vec u+f\vec v^a\\
         MO&[\langle \iota^{f\nu U}_U f g^if^{-1}\rangle, t,\vec u+f\vec v]\\
         \Sigma^\infty X&[fy,\vec u+f\vec v]\\
         B^\infty X
         & [[\langle f\ltimes \alpha^a\rangle,\langle f x^e\rangle,\vec u+f\vec v]_T,t]
         \\
         \widetilde\Omega Y&[f\vec v\mapsto \prescript{\vec u_{f\nu V^\bot}}{f}{\gamma(f^{-1}\vec u_{f\nu V}+\vec v)}]\\\hline
    \end{array}
$$

The $A$-indexed smash product is
\begin{equation*}
    \wedge_{\mathscr L A }:\Pi_A\mathds L[\texttt{Sp}]\rightarrow \mathds L[\texttt{Sp}];
    \ \wedge_{\mathscr L A}\langle Y^a\rangle :=
    \left\langle{\mathscr L A\ltimes \overline\wedge_A Y^a}_U/_{\otimes  ^{\vec u}_{f}[\prescript{\vec v^a}{g^a}{y^a} ]= \otimes^{\vec u+f_a\vec v^a}_{\langle f_ag^a\rangle}[y^a]}\right\rangle
\end{equation*}
with structural maps induced by the ones for the twisted smash product. In order to make explicit the parallel between the smash product of spectra with the tensor product of abelian groups we will use the notation
$$
    \otimes^{\vec u}_f[y^a]
    := [ [{}_f^{\vec u},[y^a] ] ]\in \wedge_{\mathscr L A}\langle Y^a\rangle,
$$
so that the $\mathds L$ structural maps are given by the formula
\begin{gather*}
    \prescript{\vec u}{f}{\otimes  ^{\vec v}_{g}[y^a]}
    :=\otimes^{\vec u+f\vec v}_{\langle fg_a\rangle}[y^a].
\end{gather*}

For $A=\underline 2$ this defines an associative and symmetric smash product 
$$
    Y^1\wedge_{\mathscr L}Y^2:=\wedge_{\mathscr L\underline 2}\langle Y^1,Y^2\rangle.
$$
Associativity follows from the fact that the maps
\begin{gather*}
    \wedge_{\mathscr L A }(\wedge_{\mathscr L B^a }\langle Y^{ab}\rangle)\rightarrow \wedge_{\mathscr L \Sigma_A B^a }\langle Y^{ab}\rangle;\ \otimes^{\vec u}_f[\otimes^{\vec v^a}_{g^a}[y^{ab}]]\mapsto \otimes^{\vec u+f_a\vec v^a}_{\langle f_ag^a_b\rangle}[y^{ab}] 
\end{gather*}
are isomorphisms \cite[Theorems I.5.4, I.5.5 and I.5.6]{EKMM}. In particular when the $B^a$ are a constant set $B$ we have a natural isomorphism
$$
    \Phi_{A,B}:\wedge_{\mathscr L A}(\wedge_{\mathscr L B}\langle Y^{ab}\rangle)\rightarrow\wedge_{\mathscr L B}(\wedge_{\mathscr L A}\langle Y^{ab}\rangle)
$$
and we set the notation
$$
    \otimes^{\vec v}_g[\otimes^{\vec u^b}_{f^b}[y^{ab}]]:=\Phi_{A,B}\otimes^{\vec u}_f[\otimes^{\vec v^a}_{g^a}[y^{ab}]].
$$

Symmetry is given by the natural isomorphism
\begin{gather*}
    \tau_{Y^1,Y^2}:Y^1\wedge_{\mathscr L} Y^2\xrightarrow\cong Y^2\wedge_{\mathscr L} Y^1,
    \hspace{0.5cm}
    \tau\otimes^{\vec u}_{\langle f_1,f_2\rangle}[y^1,y^2]
    :=\otimes^{\vec u}_{\langle f_2,f_1\rangle}[y^2,y^1] .
\end{gather*}

For all $Z\in\mathds L[\texttt{Sp}]$ set the notation $\Sigma^Z:=-\wedge_{\mathscr L}Z:\mathds L[\texttt{Sp}]\rightarrow \mathds L[\texttt{Sp}]$.

This smash product almost has a unit given by the sphere spectrum $\mathds S$, in that there are natural weak equivalences
\begin{gather*}
    \rho_Y :\Sigma^{\mathds S}Y\xrightarrow{\sim} Y,
    \hspace{0.5cm}
    \otimes_f^{\vec u}[y,\vec v]
    \mapsto\sigma^{f_1\nu U^1}_U[\prescript{\vec 0}{f_1}{y}
    ,\vec u+f_2\vec v].
\end{gather*}
Unfortunately $\rho$ is not in general a natural isomorphism, though it is a natural weak equivalence. The category of $\mathds S$-modules is the full subcategory
$$
\texttt{Mod}_{\mathds S}:=\{Y\in\mathds L[\texttt{Sp}]\mid \rho_Y\text{ is an isomorphism}\}.
$$

With the same smash product and unit $\mathds S$ the category $\texttt{Mod}_{\mathds S}$ is a symmetric monoidal category.

From the nontrivial fact that $\mathscr L A/_{\mathscr L\underline 1^A}$ has a single equivalence class \cite[Theorems I.8.1 and Section XI.2]{EKMM} the sphere spectrum $\mathds S$, the Eilenberg-Maclane spectra $HG$ and the Thom spectrum $MO$ in \ref{Example spectra}, as well as $\Sigma^{\mathds S} Y$ for $Y\in\mathds L[\texttt{Sp}]$ and spetrifications $\widetilde \Omega Y$ for $Y\in\texttt{Mod}_{\mathds S}$, are all $\mathds S$-modules with inverse maps given respectively by:
$$
    \begin{array}{|c|c|}\hline
         & \rho^{-1}y \\\hline
         \mathds S&\otimes^{\vec 0}_f[f_1^{-1}\vec u,\vec 0]\\
         HG&\otimes^{\vec 0}_f[g_a\otimes f_1^{-1}\vec u^a,\vec 0]\\
         MO&[\langle f^{-1}_1g^if_1\rangle, t,f_1^{-1}\vec u],\vec 0]\\
         \Sigma^{\mathds S}Y
         &\otimes^{\vec u}_{\langle f_1,g_2,g_3 \rangle}[y,g_2^{-1}f_2\vec v,\vec 0]
         \\
         \widetilde\Omega Y&[\vec v\mapsto \rho^{-1}\gamma\vec v]\\\hline
    \end{array}
$$
where in the first three lines $f\in\mathscr L\underline 2$ is any linear isometry such that $U\subset f_1\mathds R^\infty$, in the fourth $\langle g_2,g_3\rangle\in\mathscr L\underline 2$ are chosen such that $\langle f_1,g_2,g_3\rangle\in\mathscr L\underline 3$ and $f_2\nu U^2\subset g_2\mathds R^\infty$.

The functor $\Sigma^{\mathds S}:=-\wedge_{\mathscr L} \mathds S$ is the right adjoint of the inclusion of $\texttt{Mod}_{\mathds S}$ in $\mathds L[\texttt{Sp}]$. The functor $\Sigma^{\mathds S}$ is also a left adjoint, with right adjoint induced by a closed structure on $\mathds L[\texttt{Sp}]$ given by an $\mathds L$-mappings functor $F_{\mathscr L}$. Details of this construction can be found in \cite[Section I.7]{EKMM}, but we give an overview to establish notation. The twisted half-smash product $\mathscr L A\ltimes-$ admits a right adjoint, the \textit{twisted function spectrum} functor

\begin{gather*}
    F[\mathscr L A ,-):\texttt{Sp}\rightarrow \texttt{Sp}_{\oplus_A\mathds R^\infty},
    F[\mathscr L(A),Y):=
    \left\langle \underset{K\subset_{cpct}\mathscr L(A)}{\text{Lim }} Y_{\mu  \langle U^a\rangle}^{TK_{\mu \langle U^a\rangle}^{\langle U^a\rangle}}\right\rangle,\\
    \sigma^{\langle V^a\rangle}_{\langle U^a\rangle}[\varphi,\langle \vec v^a\rangle]
    :=\left\langle\prescript{\vec u}{f}{}\mapsto\sigma^{\mu\langle V^a\rangle}_{\mu\langle U^a\rangle}\left[\varphi_f^{(\vec u+f_a\vec v^a)_{\mu \langle U^a\rangle}},(\vec u+f_a\vec v^a)_{\mu\langle U^a\rangle^\bot}\right]\right\rangle.
\end{gather*}

For $U^1\in\mathscr A$ we also have a shift functor
\begin{gather*}
    -[U^1]:\texttt{Sp}_{\mathds R^\infty\oplus\mathds R^\infty}\rightarrow \texttt{Sp}; Y[U^1]=\langle Y_{U^1,U^2}\rangle, \sigma^{U^2}_{V^2}[y,\vec v]:=\sigma^{U^1,U^2}_{U^1,V^2}[y,(\vec 0,\vec v)].
\end{gather*}

If $Y\in\mathds L[\texttt{Sp}]$ then $F[\mathscr L\underline 2,Y)[U^1]\in \mathds L[\texttt{Sp}]$ with structural map
$$
    \prescript{\vec u}{f}{\varphi}:=
    \langle\prescript{\vec v}{g}{}\mapsto
    \varphi^{\vec v+g_2\vec u}_{\langle g_1,g_2f\rangle}
    \rangle.
$$

Finally, we can now define
\begin{gather*}
    F_{\mathscr L}(-,-):\mathds L[\texttt{Sp}]^{\text{op}}\times\mathds L[\texttt{Sp}]\rightarrow\mathds L[\texttt{Sp}];\\
    F_{\mathscr L}(Z,Y):=
    \left\langle\left\{\phi\in\mathds L[\texttt{Sp}](Z, F[\mathscr L \underline 2 ,Y)[U^1])\left\lvert
    \begin{array}{l}
         \prescript{\vec u}{f}{(\phi z^{\vec v}_{g})}=\phi z^{\vec u+f\vec v}_{\langle fg_1,fg_2\rangle},\\
         \phi(\prescript{\vec u}{f}{z})^{\vec v}_g
         = \phi z^{\vec u+g_2\vec v}_{\langle g_1,g_2f\rangle}.
    \end{array}\right.\right\}\right\rangle,\\
    \sigma^{U^1}_{V^1}[\phi,\vec v]:=
    \langle z\mapsto \sigma^{U^1,U^2}_{V^1,U^2}[\phi z,(\vec v,\vec 0)]\rangle,
    \hspace{0.5cm}\prescript{\vec u}{f}{\phi}:=\langle [z,\prescript{\vec v}{g}{}]\mapsto \phi z^{\vec v+g_1\vec u}_{\langle g_1f,g_2\rangle}\rangle.
\end{gather*}

The functor $F^{\mathds S}:=F_{\mathscr L}(\mathds S,-):\texttt{Mod}_{\mathds S}\rightarrow \mathds L[\texttt{Sp}]$ is right adjoint to $\Sigma^{\mathds S}$. 

The monoidal structure of $\mathds S$-modules provides a natural definition of ring spectra, module spectra and algebra spectra.

\begin{definition}
    A \textit{commutative ring spectrum} $R$ is a commutative monoid in $\texttt{Mod}_{\mathds S}$, ie an $\mathds S$-module equipped with a \textit{multiplication} map $\mu:R\wedge_{\mathscr L}R\rightarrow R$ and a \textit{unit} map $\eta:\mathds S\rightarrow R$ satisfying natural associative, unit and commutative laws. The category of commutative ring spectra is denoted $\texttt{CRingSp}$. 
    
    For $R\in\texttt{CRingSp}$ an \textit{$R$-module} $M$ is a module over $R$, ie an $\mathds S$-module equipped with an \textit{action} $\lambda : R\wedge_{\mathscr L} M\rightarrow M$, satisfying natural associative and unit laws.  The category of $R$-modules is denoted as $\texttt{Mod}_R$. 
    
    The category of $R$-modules admits a monoidal structure with associative and symmetric tensor product the coequalizer
    $$
    M\wedge_RN:=
    \text{Coeq}(M\wedge_{\mathscr L}R\wedge_{\mathscr L}N\rightrightarrows M\wedge_{\mathscr L} N)
    $$
    and unit $R$. The category of $R$-modules is denoted $\texttt{Mod}_R$.
    
    A commutative $R$-algebra is a commutative monoid in $(\texttt{Mod}_R,\wedge_R,R)$, and the category of commutative $R$-algebra is denoted $\texttt{CAlg}_R$.
    
    The category of commutative algebra spectra is defined as 
    $$\texttt{CAlgSp}:=\Sigma_{\texttt{CRingSp}}\texttt{CAlg}_R.
    $$
\end{definition}

As in the classical set theoretical setting there is a natural isomorphism
$\texttt{CAlgSp}\cong \texttt{CRingSp}^\shortrightarrow$ \cite[VII.1]{EKMM}. Alternatively we have a monad $(\mathds P^\shortrightarrow;\eta,\mu)$ on $\mathds L[\texttt{Sp}]^2$ with
\begin{gather*}
    \textstyle \mathds P^\shortrightarrow Y_\star:=\int^{\mathds S_\star}\wedge_{\mathscr L A} \langle Y_{\mathfrak c a}\rangle;\\ 
    \eta_\star y:=[\otimes^{\vec 0}_{id} y ],
    \hspace{1cm} 
    \mu_\star[\otimes^{\vec u}_{f}[\otimes^{\vec v^a}_{g^a}[y^{a,b}]]]:=[\otimes^{\vec u+f_a\vec v^a}_{\langle f_ag^a_b\rangle}[y^{a,b}]]
\end{gather*}
which restricts to a monad on $\texttt{Mod}_{\mathds S}^2$. The objects of $\mathds P^\shortrightarrow[\mathds L[\texttt{Sp}]^2]$ behave like algebra spectra over ring spectra except they have units only up to weak equivalence and are refered to as $E^\shortrightarrow_\infty$-algebra spectra, similarly to how algebras in $\mathds L[\texttt{Sp}]$ over the nonrelative version $\mathds P$ of this monad are called $E_\infty$-ring spectra. By the same argument as in \cite[Prop. II.4.5]{EKMM} we have an isomorphism 
\begin{equation}\label{CAlgSpArePAlg}
    \mathds P^\shortrightarrow[\texttt{Mod}_{\mathds S}^2]\cong\texttt{CAlgSp}.
\end{equation}
For $R=((R_d,R_c);\eta,\mu)\in \texttt{CAlgSp}$ and $\otimes^{\vec u}_f[r^a]\in \wedge_{\mathscr LA}\langle R_{\mathfrak c a}\rangle$ we will use the notation
$$
    \textstyle\prod^{\vec u}_f[r^a]:=\mu[\otimes^{\vec u}_f[r^a]].
$$

The sphere spectrum $\mathds S$, the Eilenberg-MacLane spectrum of a commutative ring $HR$, the Thom spectrum $MO$, suspensions $\Sigma^\infty X$ of $\mathscr L$-spaces, deloopings $B^\infty X$ of $E_\infty$-rings, the $\mathds S$-module $\Sigma^{\mathds S} R$ associated to an $E_\infty$-ring spectrum $R$ in $\mathds P[\mathds L[\texttt{Sp}]]$ and spectrifications $\widetilde \Omega R$ of ring spectra $R$ are all commutative ring spectra with
$$
    \begin{array}{|c|c|c|}\hline
         & \eta\vec u& \prod^{\vec u}_f[y^a]\\\hline
         \mathds S&\vec u&
         \vec u+f_a\vec v^a\\\hline
         HR&1_R\otimes \vec u&\prod_{A_{\langle b^a\rangle}} r^a_{b^a}\otimes \vec u+f_a\vec v^{a,b^a}\\\hline
         MO&[id,\emptyset,\vec u]&[\prod_A\langle \iota^{f_a\nu U^a}_U  f_ag^{a,i}f_a^{-1}\rangle\cdot\delta^a,\lhd_A t^a,\vec u+f_a\vec v^a]\\\hline
         \Sigma^\infty X&[1_X,\vec u]&[f[x^a],\vec u+f_a\vec v^a]\\\hline
         B^\infty X&\![[\emptyset,1_X,\vec u]_{\underline 1},\emptyset]&
         \!\left[\left[\!\left\langle f_a\ltimes \left\langle
         \begin{cases}
         \textbf{id}\\
         \alpha^{a,v}
         \end{cases}\hspace{-0.4cm}\right\rangle\right\rangle,\langle f_a[x^{a,e}]\rangle,\vec u+f_a\vec v^a\!\right]_{\Pi_A T^a\cdot\delta^a}\hspace{-0.7cm},\lhd_A t^a\!\right]\\\hline
         \Sigma^{\mathds S}R&\otimes^{\vec u}_f[1_Y,\vec u]&
         \otimes^{\vec v}_{g}[\prod^{\vec u^1}_{f^1}[r^a],\vec u^2+f^2_a\vec v^a]
         \\\hline
         \widetilde\Omega R&[\vec v\mapsto \sigma^U_V[\eta\vec u,\vec v]]&[f_a\vec v^a\mapsto \prod^{\vec u_{f\nu V^\bot}}_f[\gamma^a(f_a^{-1}\vec u_{f_a\nu V^a}+\vec v^a)]]\\\hline
    \end{array}
$$
with the implicit conditions in the fifth line as in the formula \ref{LActB2X}.

There is a natural isomorphism $\texttt{CAlg}_{\mathds S}\cong\texttt{CRingSp}$, which is analogous to the isomorphism between commutative rings and commutative $\mathds Z$-algebras. Moreover $(MO,HR)\in \texttt{CAlgSp}$ with
$$
\textstyle\prod^{\vec u}_f[
[\langle g^{1,i}\rangle, t^1, \vec v^1],
r^2_b\otimes \vec v^{2,b}]
:= r^2_b\otimes \vec u+ f_1\circ_{\underline m}g^{1,i}\vec v^1 + f_2\vec v^{2,b}.
$$

\subsection{Stable mixed model structure of commutative algebra spectra}

The stable mixed model structure of $\texttt{Mod}_{\mathds S}$ is right transferred from the one in $\texttt{Sp}$ by the adjunction
$$
(\Sigma^{\mathds S}\mathds L\dashv F^{\mathds S}):\texttt{Sp}\leftrightharpoons \texttt{Mod}_{\mathds S}
$$
as described in \cite{BR13,Cole06,EKMM}, so that weak equivalences and fibrations in $\texttt{Mod}_{\mathds S}$ are those maps whose underlying spectrum mapping are $q$-equivalences and $h$-fibrations respectively. The Hurewicz/Str{\o}m factorization systems are constructed as in \texttt{Sp} with the $\mathds S$-module structures of $\Gamma \mathfrak f$ and $E\mathfrak f$ defined point-wise. The mixed model structure of $\texttt{CAlgSp}$ is right transferred from the one in $\texttt{Mod}^2_{\mathds S}$ by the adjunction
$$
(\mathds P^\shortrightarrow\dashv U): \texttt{Mod}_{\mathds S}^2\leftrightharpoons\texttt{CAlgSp}.
$$
The Quillen model structure is transferred due to the fact that $\texttt{CAlgSp}$ has continuous coequalizers and satisfies the “Cofibration
Hypothesis” as in \cite[Theorem VII.4.7]{EKMM}. The Hurewicz/Str\o m model structure is transferred since we can define an algebra structure on $\Gamma \mathfrak f$ for $\mathfrak f\in \texttt{CAlgSp}$ as
$$
\textstyle \eta\vec u:=(\eta_X\vec u,0,r\mapsto \eta_Y\vec u);\hspace{0.5cm}
\prod^{\vec u}_f[(x^a,t^a,\gamma^a)]:=
(\prod^{\vec u}_f x^a,\max_A t^a,r\mapsto \prod^{\vec u}_f\gamma^a r)
$$
As in $\texttt{Sp}$ we have that $(\Gamma; C_t,F)$ forms an algebraic weak factorization system in $\texttt{CAlgSp}$. On the other hand there doesn't seem to be any natural algebra structure on $E\mathfrak f$ such that the $h$-cofibration/trivial $h$-fibration factorization $(E;C,F_t)$ in $\texttt{Sp}$ induces a factorization in $\texttt{CAlgSp}$. We do have an $h$-cofibration/$h$-equivalence factorization
$$
\xymatrix@C=2cm{X\ar@{^{(}->}[r]^(0.25){\text{in}_X}&X\wedge_{\mathds P X}\mathds P(\Gamma\mathfrak f\wedge [0,\infty]_+)\wedge_{\mathds P\Gamma\mathfrak f} Y
\ar[r]_(0.75)\cong^(0.75){(\mathfrak f,F_t\mathfrak f^\dagger,id)}&Y}
$$
and the fact that $C_t\mathfrak f$ has the left lifting property against $h$-fibrations in $\texttt{Mod}_{\mathds S}$ induces the left lifting property against $h$-fibrations in $\texttt{CAlgSp}$ on $\text{in}_X$. The map $(\mathfrak f,F_t\mathfrak f^\dagger,id)$ is an $h$-equivalence, but it is not necessarily an $h$-fibration. Applying $(\Gamma; C_t, F)$ then gives us the $h$-cofibration/trivial $h$-fibration factorization
$$\xymatrix@C=2.7cm{
    X\ar@{^{(}->}[r]^(0.4){C_t(\mathfrak f,F_t\mathfrak f^\dagger,id)\text{in}_X}&
    \Gamma (\mathfrak f,F_t\mathfrak f^\dagger,id)\ar@{->>}[r]_(0.6)\cong^(0.6){F(\mathfrak f,F_t\mathfrak f^\dagger,id)}&Y
}$$
which determines the Hurewicz/Str{\o}m, and therefore also the mixed, model structure on $\texttt{CAlgSp}$.

\subsection{Recognition of algebra spectra}

Let $\mathcal E^\shortrightarrow$ be an $E^\shortrightarrow_\infty$-operad equipped with an $\mathscr L^\shortrightarrow$-action. The functors $F^{\mathds S}$ and $\Sigma^{\mathds S}$ induces objectwise adjoint functors $F^{\mathds S}_\shortrightarrow$ and $\Sigma^{\mathds S}_\shortrightarrow$ on the morphism categories. We can then define the functors
\begin{gather*}
    \Omega_2^{\infty,\mathds S}:\texttt{CAlgSp}\rightarrow
    (\mathscr L^\shortrightarrow,\mathcal E^\shortrightarrow)[\texttt{Top}], 
    \hspace{0.5cm}\Omega_2^{\infty,\mathds S}R:= \Omega^\infty_2 F^\mathds S_\shortrightarrow \eta;\\
    f[\phi^a]:=[\vec u^1\mapsto\langle [\vec u^2,\prescript{\vec v}{g}{}]\mapsto
    \textstyle\prod_f^{\vec 0}\phi^a[\vec u^1,\vec u^2]^{\vec v}_g
    \rangle]
\end{gather*}
and
\begin{gather*}
    B^{\mathds S,\infty}_\shortrightarrow:(\mathscr L^\shortrightarrow,\mathcal E^\shortrightarrow)[\texttt{Top}]\rightarrow\ \texttt{AlgSp};
    \hspace{0.5cm} B_\shortrightarrow^{\mathds S,\infty}X:=\Sigma^{\mathds S}_\shortrightarrow B^\infty_\shortrightarrow X.
\end{gather*}

\begin{theorem}\label{RecogCAlgSp}
    There is an idempotent quasiadjunction 
    $$
    (B^{\mathds S,\infty}_\shortrightarrow\dashv_{\ \overline B_2,\widetilde\Omega_\shortrightarrow } \Omega^{\infty,\mathds S}_2):(\mathscr H^\shortrightarrow_\infty,\mathscr L^\shortrightarrow)[\texttt{Top}]\leftrightharpoons\texttt{CAlgSp}
    $$
    that induces an equivalence of homotopy categories
    $$
        (\mathbb LB^{\mathds S,\infty}_\shortrightarrow\dashv \mathbb R\Omega^{\infty,\mathds S}_2):\mathcal Ho (\mathscr H^\shortrightarrow_\infty,\mathscr L^\shortrightarrow)[\texttt{Top}]_{\text{alg}}\leftrightharpoons\mathcal Ho \texttt{CAlgSp}_{\text{con}}.
    $$
\end{theorem}
  
\textbf{Proof:} The natural weak equivalences $\eta'$ and $\epsilon'$ are defined as in the proof of theorem \ref{MorphRecogPrin}. The other natural transformations of the unit span and counit cospan are
\begin{gather*}
    \eta_\star\![[\!\alpha^r\!,\!\langle\alpha^v\rangle\!,\!\langle x^e\rangle\!]_T\!,\! t\!]
    \!:=\!\left[\!\vec u^1\!\!\mapsto\!\begin{cases}
        \langle [\vec u^2,\prescript{\vec v}{f}{}]\mapsto\otimes_f^{\vec v}[[[\langle\alpha^v\rangle,\langle x^e\rangle,\alpha^{r,-1}_{e'}\vec u^1]_{T_{\geq e'}}, t],\vec u^2]\rangle\\
        \infty
    \end{cases}\!\!\!\!\!\!\!\right];\\
    \epsilon_{\star,U}
    \otimes^{\vec u}_f
    [[[
    \langle\alpha^v\rangle,
    \langle \phi^e\rangle,\vec v^1\!]_{T},
     t],\vec v^2\!]
    \!:=\! [\vec w\!\mapsto\!
    \circ_T\alpha^v\langle \phi^e\rangle[\vec v^1+\vec w^{f_1},
    \vec v^2+\vec w^{f_2}]_f^{\vec u+\vec w_{f^\bot}}\!].
\end{gather*}
with the conditions in the first formula as in the proof of theorem \ref{MorphRecogPrin} and with the domain in the last formula any $W\in\mathscr A_{U+f_1V}$ with $V$ a common domain of representatives of the loops $\phi^e$.

That these maps satisfy the conditions for an idempotent quasiadjunction follows from the fact that $(\Sigma^{\mathds S}_\shortrightarrow\dashv F^{\mathds S}_\shortrightarrow):\mathds P^\shortrightarrow[\mathds L[\texttt{Sp}]^2]\leftrightharpoons\texttt{CAlgSp}$ is a Quillen equivalence and the same argument as for \ref{MorphRecogPrin} and \ref{Idempotent}.$\blacksquare$

\bibliographystyle{amsplain}
\bibliography{bibliography.bib}
\end{document}